\scrollmode
\documentclass[12pt]{amsart}

\usepackage{amssymb}
\usepackage{eepic}
\usepackage{mathrsfs}
\usepackage{xypic}
\usepackage{verbatim}
\usepackage{graphicx}
\usepackage{psfrag}
\CompileMatrices

\def\makemargins{
	\oddsidemargin .25in
	\evensidemargin .25in
	\textwidth 6.00in
}
\makemargins
\theoremstyle{plain}
\newtheorem{theorem}[subsection]{Theorem}
\newtheorem{lemma}[subsection]{Lemma}
\newtheorem{proposition}[subsection]{Proposition}

\theoremstyle{definition}
\newtheorem{definition}[subsection]{Definition}

\theoremstyle{remark}
\newtheorem{remark}[subsection]{Remark}

\newcommand{\R}{{\mathbb R}}
\newcommand{\C}{{\mathbb C}}
\newcommand{\bbA}{{\mathbb A}}

\newcommand{\cE}{{\mathcal E}}
\renewcommand{\cH}{{\mathcal H}}
\newcommand{\sO}{{\mathscr O}}
\newcommand{\sF}{{\mathscr F}}
\newcommand{\sB}{{\mathscr B}}
\newcommand{\sC}{{\mathscr C}}

\newcommand{\core}{{Z}}
\newcommand{\proj}{{\mathbb P}}

\renewcommand{\tilde}{\widetilde}

\newcommand{\card}[1]{|#1|}
\newcommand{\closure}[1]{\overline{#1}}
\renewcommand{\emptyset}{\varnothing}

\newcommand{\Grass}{\operatorname{Gr}}
\newcommand{\Flag}{\operatorname{Fl}}
\newcommand{\Dim}{\operatorname{dim}}

\newcommand{\Proj}{\operatorname{Proj}}
\newcommand{\Spec}{\operatorname{Spec}}
\newcommand{\Conv}{\operatorname{Conv}}
\newcommand{\SL}{\operatorname{SL}}
\newcommand{\GL}{\operatorname{GL}}
\newcommand{\mainset}{[\![4]\!]}

\catcode`\@=11
\def\@secnumfont{\bfseries}
\catcode`\@12

\begin{document}
\title{A smooth space of tetrahedra}

\newif \ifdraft

\def \makeauthor{
\author{Eric Babson}
\address{Department of Mathematics\\
University of Washington\\
Seattle, WA 98195}
\email{babson@math.washington.edu}

\author{Paul E. Gunnells}
\address{Department of Mathematics\\
Columbia University\\
New York, NY 10027}
\email{gunnells@math.columbia.edu}

\author{Richard Scott}
\address{Department of Mathematics\\
Santa Clara University\\
Santa Clara, CA 95053}
\email{rscott@math.scu.edu}
}
\makeauthor

\date{October 1999.  Revised June 2000}
\subjclass{14M15}
\keywords{Compactifications of configuration spaces, space of tetrahedra, space
of triangles.}

%%%
%%% Abstract
%%%
\begin{abstract}
We construct a smooth symmetric compactification of the space of all labeled
tetrahedra in $\proj ^{3}$.
\end{abstract}

\maketitle

%%%
%%% Introduction
%%%
\section{Introduction}

\subsection{}
Let $\proj ^{n} $ be $n$-dimensional complex projective space, and let
$P\subset \proj ^{n}$ be a set of $n+1$ labeled points in general
position.  By taking all possible linear spans of subsets of $P$,
one obtains a configuration of flats in $\proj ^{n}$ arranged to form
a simplex.  The set $X^{\circ }$ of all such configurations is
naturally a quasi-projective variety with a canonical singular
compactification $X$.  One is interested in the variety $X$ for many
reasons: 

\begin{itemize}
\item For $n=2$, the space $X$ is the space of triangles in the
plane.  In \cite{schubert}, Schubert described a
desingularization of $X $, and used it to study enumerative problems
involving triangles \cite{semple, speis-robts, collino-fulton}.  
\item The space $X $ is a \emph{configuration variety} in the sense
of Magyar \cite{magyar1, magyar2}.  Such spaces arise naturally in the
study of \emph{generalized Schur modules}.  These are (reducible)
$\GL_{n}$-modules that generalize the classical Schur modules, and have
been studied in various guises by many authors \cite{buchsbaum, kras2,
woodcock, kras1, vic2, vic1, shi}.  One hopes that configuration varieties
will play a role in a Borel-Weil theory for these modules.
\item Let $\sB$ be the Tits building for $\SL_{n+1} (\C
)$, and let $\sC$ be the associated Coxeter complex \cite{tits}.
Then $X $ can be interpreted as the space of maps of $\sC$ into
$\sB$.  By considering other algebraic groups, one obtains a
collection of natural configuration spaces related to the \emph{Bott-Samelson
varieties} of Demazure \cite{demazure}.  In particular, $X $ can be
regarded as a \emph{canonical} Bott-Samelson variety associated to
all reduced expressions of the longest word of the Weyl group of
$\SL_{n+1}$.
\item The space $X $ is a natural generalization of the
Fulton-MacPherson space \cite{fulton-macpherson} $\proj
^{n}[n+1]$.  This variety adds data to an open set of the product
$\prod _{i=1}^{n+1}\proj ^{n}$ that records how points approach the
(large) diagonals.  In fact, $\proj^n[n+1]$ is a desingularization of
the space of all \emph{$1$-skeleta} of $n$-simplices in $\proj
^{n}$.
\end{itemize}

\subsection{}
\begin{figure}[ht]
\begin{center}
\includegraphics[scale = .5]{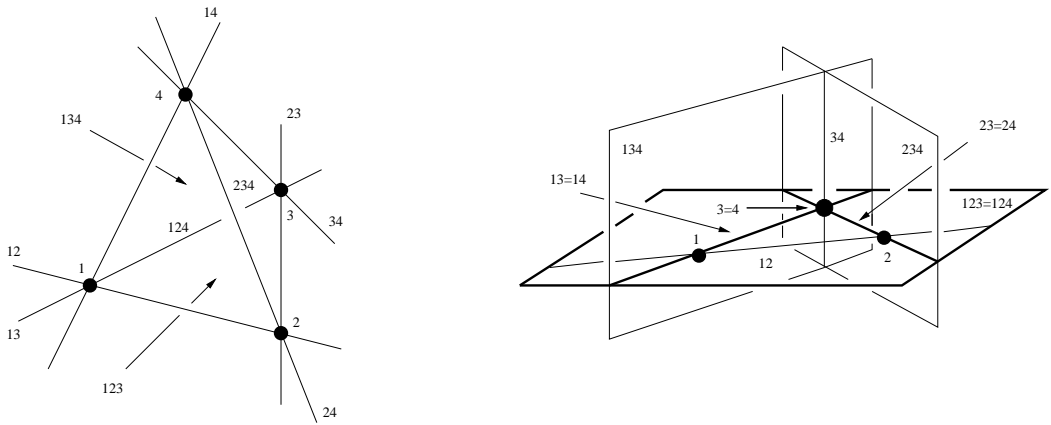}
\end{center}
\caption{\label{tet.fig} A point in $X^{\circ}$ and a point in $X\smallsetminus X^{\circ }$}
\end{figure}

In this paper we consider the case $n=3$, where $X$ is the space of
tetrahedra in $\proj ^{3}$.  We construct a symmetric compactification
$\tilde{X}$ of $X^{\circ} $ that we call the space of \emph{complete
tetrahedra}.  It is obtained by embedding $X^{\circ}$ into a large
ambient variety $E_{\#}$ and taking the closure of the image.  The
singular locus of the canonical compactification $X$ is contained in
the subvariety consisting of ``collapsed tetrahedra''---that is,
configurations of flats where certain faces coincide
(Figure~\ref{tet.fig})---and $E_{\#}$ is constructed to capture the
asymptotic behavior of a tetrahedron as it collapses.  Our main
theorem (Theorem \ref{mainthm}) is that $\tilde{X}$ is nonsingular.

In a later paper \cite{tet.coho} we will
make a more detailed study of the geometry of $\tilde{X}$.  We will
show that the complement of $X^{\circ }$ in $\tilde{X}$
is a divisor with normal crossings and compute the cohomology ring
of $\tilde{X}$.  

\subsection{}

Although this article considers the space of tetrahedra in $\proj
^{3}$, the definition of $\tilde{X}$ makes sense for all $n$.
Many of the results of the paper (in particular,
sections~\ref{s:infinity}--\ref{s:core}) hold for arbitrary $n$, but
we have avoided this generality since we cannot complete the proof
that $\tilde{X}$ is nonsingular in general (the combinatorial arguments in
section~\ref{s:nonsing} become infeasible when $n\geq 4$).  However, we
conjecture that $\tilde{X}$ provides a nonsingular compactification
of $X^{\circ }$ for all $n$.

For $n=2$, it is not hard to see that our variety $\tilde{X}$ is
nonsingular and coincides with certain triangle varieties found in the
literature.  More precisely, it is isomorphic to the Fulton-MacPherson
space $\proj ^{2}[3]$, which in turn coincides with an auxiliary
compactification constructed by Roberts-Speiser \cite{speis-robts}.
It is not, however, isomorphic to Schubert's compactification as a
variety over $X$.\footnote{The difference between Schubert's space and
Fulton-MacPherson-Roberts-Speiser's space appears when one considers
the torus action on them, cf. \cite{magyar2}}

\subsection{} 
We now give an overview of the definition of $\tilde{X}$.  The
construction of $E_{\#}$ depends on the combinatorics of
\emph{hypersimplices} \cite{gelfand-macpherson}, polytopes intimately
related to the geometry of Grassmannians.  For our considerations, the
relevant polytopes are the $3$-dimensional hypersimplices $\Delta
_{1}$, $\Delta _{2}$, and $\Delta _{3}$ (Figure~\ref{labgeneric.fig}).
The vertices of these hypersimplices are in bijection with the labeled
faces of a tetrahedron, and the edges of the hypersimplices correspond
to certain pairs of faces of the same dimension.

For each edge $\alpha$ in a hypersimplex, we form a $(\proj ^{1}\times
\proj^{1})$-bundle $E_{\alpha }\rightarrow X$.  The bundle $E_{\alpha
}$ has a canonical section $u_{\alpha }$ and a diagonal subbundle
$D_{\alpha}$.  The section $u_{\alpha}$ tracks the subspaces
corresponding to the vertices of $\alpha$ in such a way that
$u_{\alpha}$ intersects $D_{\alpha}$ precisely when these subspaces
coincide.  In order to record the asymptotic behavior in $X$ near a
collapsed tetrahedron, a natural idea is to form products of the
$E_{\alpha}$'s and blowup the corresponding product of diagonals.  The
question is which products to take, and why.

\subsection{}
Our main idea is that the relevant products are those indexed by the
faces of dimension $\geq 2$ of the hypersimplices.  The motivation is
that a configuration of flats in $\proj ^{3}$ arranged to form a
tetrahedron contains certain ``sub-'' and ``quotient'' configurations
corresponding to proper faces of the hypersimplices.  For example, the
three points and three lines in a given face of a tetrahedron form a
subconfiguration that corresponds to a triangular face in the
hypersimplex $\Delta_1$, and the three lines and three planes
containing a given point form a quotient configuration that
corresponds to a triangular face in the hypersimplex $\Delta_3$.  Our
motivation is that a nonsingular compactification of $X^{\circ }$
should add data recording the ``infinitesimal shapes'' of these sub-
and quotient configurations.  Hence, each locus we blowup corresponds
to the collapsing together of the subspaces labeled by some face of a
hypersimplex.

More precisely, our definition is as follows.  Let $\cH $ be the set
of faces of dimension $\geq 2$ of all the $\Delta _{k}$.  For each
$\beta \in \cH $, let $\cE (\beta )$ be the set of edges in $\beta $.
Let $E _{\beta }$ be the product bundle
\[ 
E _{\beta } := \prod _{\alpha \in \cE (\beta) } E _{\alpha },
\]
and let $D_{\beta}$ be the corresponding product of the diagonals
$D_{\alpha}$.  The ambient variety $E_{\#}$ is then defined by blowing
up each $E_{\beta}$ along $D_{\beta}$ and taking the product of the
resulting blowups.  The corresponding product of the sections
$u_{\alpha}$, then determines an embedding $X^{\circ }\rightarrow
E_{\#}$, and we define $\tilde{X}$ to be the closure.

\subsection{}
The paper is organized as follows.  Section~\ref{s:setup} sets up
notation, defines $X$, and contains background on hypersimplices.
Section~\ref{s:global} contains the construction of $E_{\#}$ and
$\tilde{X}$.  In section~\ref{s:infinity} we describe a collection of
affine open subvarieties that covers $X$ and give equations defining a
typical element $U\subset X$ in this collection.  In
section~\ref{s:infinity-tilde} we restrict the construction of
$E_{\#}$ to $U$ to obtain $\tilde{U}$, a certain subvariety of
$\tilde{X}$.  The point of
sections~\ref{s:infinity}--\ref{s:infinity-tilde} is that the
nonsingularity of $\tilde{X}$ follows from the nonsingularity of
$\tilde{U}$.

In the remaining sections we prove nonsingularity of $\tilde{U}$.
First, in section~\ref{s:core}, we show that $\tilde{U}$ has the
structure of a vector bundle over a certain (multi-) projective
variety $\core $ that we call the \emph{core}.  We then study the
$\GL_{4}$-action on $\tilde{U}$ to show that nonsingularity of $\core
$ follows from its nonsingularity at points in a certain subvariety
$\core _{sp}\subset \core $.  Finally in section~\ref{s:nonsing}, we
give equations that cut out $\core $ from projective space and use a
graphical description of these defining relations to show that
$\core_{sp}$ consists of nonsingular points of $\core$; this proves
Theorem~\ref{mainthm}.

\subsection{Acknowledgments}
We thank Michael Thaddeus for suggesting the proof of
Lemma~\ref{3flat-lemma}.  We thank Robert MacPherson, who originally
told us about this problem, in terms of buildings and Coxeter
complexes, and who has offered us much encouragement and interest.
Finally, we thank the various institutions that have supported us and
have hosted our collaboration at one time or another: Columbia
University, Cornell University, the Institute for Advanced Study, the
Ohio State University, and the University of Washington.

%%%
%%% Setup
%%%
\section{Notation and the basic variety $X$}\label{s:setup}

\subsection{}\label{ss:maindef} 

Let $e_1,\dots ,e_4$ be the standard basis of $\C^4$, and let
$\mainset$ be the set $\{1,2,3,4\}$.  For any subset
$I\subset\mainset$, let $E_I\subset \C ^{4}$ be the subspace spanned
by $\{e_i\mid i\in I\}$.  Let $\proj^{3}$ be the projective space of
lines in $\C^4$, and let $G$ be the algebraic group $\GL_4(\C )$.

For $k=1,2,3$, let $\Grass_{k}$ be the Grassmannian of $k$-dimensional
subspaces of $\C^4$, and for each proper nonempty subset
$I\subset\mainset$, let $\Grass_I:=\Grass_{\card{I}}$.  (We use the
notation $\card{I}$ for the cardinality of $I$.) Let $Y$ be the
product
\[
Y := \prod_{{\varnothing \subsetneq I\subsetneq\mainset}}\Grass
_{I} \cong (\Grass_1)^4\times(\Grass_2)^6 \times(\Grass_3)^4,
\] 
and for each $I\subset\mainset$, let $\pi_I$ be the projection to the
$I$th factor. The group $G$ acts on $Y$ by left multiplication, and
each $\pi _{I}$ is $G$-equivariant. 

\begin{definition}\label{the.space}
Let $p_0\in Y $ be the point such that $\pi_I(p_0)=E_I$ for all
$I\subset\mainset$, and let $X^{\circ } \subset Y $ be the $G$-orbit
of $p_{0}$.  Let $X\subset Y$ be $\closure{X^{\circ }}$ (the bar
denotes Zariski closure).  The space $X$ (respectively $X^{\circ }$)
is called the \emph{space of tetrahedra (resp. nondegenerate
tetrahedra)}.
\end{definition}

Note that since the $G$-action preserves incidence relations among
subspaces, for any $p\in X$ we have $\pi_I(p)\subset\pi_J(p)$ if
$I\subset J$.  Hence the configuration of subspaces $\{\pi_I(p)\mid
I\subset\mainset \}$ satisfies the incidence relations corresponding
to the faces of a tetrahedron.

The symmetric group $S_4$ acts on $Y$ by permuting the factors, and
this clearly induces an action on $X^{\circ}$ and $X$: given $\sigma \in
S_{4}$ and $p\in X$, the point $\sigma\cdot p$ is determined by
$\pi_I(\sigma\cdot p)=\pi_{\sigma^{-1}(I)}(p)$.  This action can be
viewed as ``changing the labels'' on the faces of a tetrahedron.

\subsection{}\label{hyper}
The construction of our resolution 
$\tilde{X}\rightarrow X$ 
is
based on the combinatorics of hypersimplices, so we recall basic facts
about them.  More details can be found in \cite{gelfand-macpherson}.

\begin{figure}[ht]
\begin{center}
\psfrag{1}{$\scriptscriptstyle{1}$}
\psfrag{2}{$\scriptscriptstyle{2}$}
\psfrag{3}{$\scriptscriptstyle{3}$}
\psfrag{4}{$\scriptscriptstyle{4}$}
\psfrag{12}{$\scriptscriptstyle{12}$}
\psfrag{13}{$\scriptscriptstyle{13}$}
\psfrag{14}{$\scriptscriptstyle{14}$}
\psfrag{23}{$\scriptscriptstyle{23}$}
\psfrag{24}{$\scriptscriptstyle{24}$}
\psfrag{34}{$\scriptscriptstyle{34}$}
\psfrag{123}{$\scriptscriptstyle{123}$}
\psfrag{124}{$\scriptscriptstyle{124}$}
\psfrag{134}{$\scriptscriptstyle{134}$}
\psfrag{234}{$\scriptscriptstyle{234}$}
\includegraphics[scale = .6]{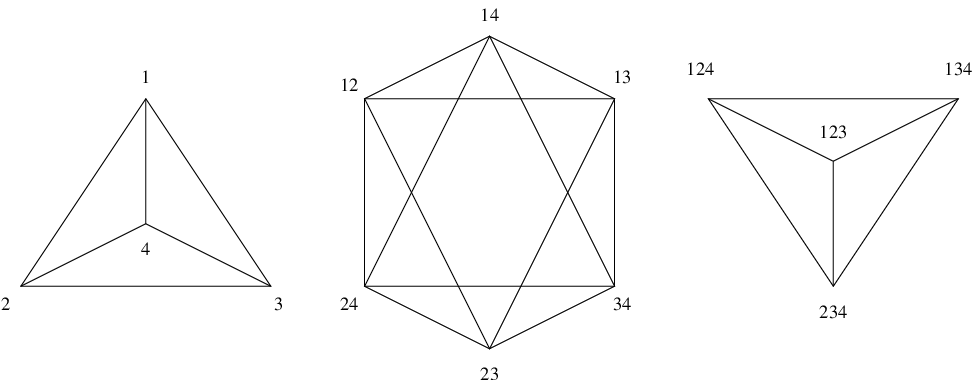}
\end{center}
\caption{\label{labgeneric.fig}}
\end{figure}

Let $\varepsilon _1,\ldots, \varepsilon _4$ be the standard basis of
$\R^{4}$, and for any subset $I\subset \mainset$, let $\varepsilon _{I} := \sum
_{i\in I} \varepsilon _{i}$.  Then the \emph{hypersimplex (of rank $k$)}
$\Delta_k$ is defined by
\[
\Delta _{k} := \Conv\bigl\{\varepsilon _{I} \bigm| I\subset
\mainset\quad \hbox{and}\quad \card{I}=k\bigr\},
\]
where $\Conv$ denotes convex hull.  The hypersimplices $\Delta_1$ and
$\Delta_3$ are $3$-simplices, and $\Delta_2$ is an octahedron.  Note
that the vertices of a hypersimplex are indexed by proper nonempty
subsets of $\mainset$ (Figure~\ref{labgeneric.fig}).  It will be
convenient to fix a total ordering on the subsets of $\mainset$, and
thus on the vertices of the hypersimplices:
\[\emptyset<1<2<3<4<12<13<14<23<24<34<123<124<134<234<1234.\]
If $I$ is a subset of $\mainset$ with $\card{I}=k$, we let $I_{0}$ be
the set $1$, $12$, or $123$, depending on whether $k=1,2$, or $3$,
respectively.

\subsection{} 
We identify faces of the hypersimplices with their corresponding sets
of vertices.  Let $\cE$ be the set of pairs $\{I,J\}$ with
$I,J\subset\mainset$ corresponding to edges of the hypersimplices,
and for $k=1,2,3$ let $\cE_k\subset \cE $ be the subset corresponding
to edges of $\Delta_k$.  Let $\cH$ be the set of vertex sets of all
faces of dimension $\geq 2$.  Hence $\cH$ contains the maximal
$3$-dimensional faces $\{1,2,3,4\}$, $\{12,13,14,23,24,34\}$, and
$\{123,124,134,234\}$, as well as $16$ triangular faces.  These last
faces can be oriented as follows.  Let $\beta = \{I,J,K \}\in \cH $ be
a triangular face, and let $(JK,IK,IJ)$ be the corresponding triple of
edges.  We call this triple an {\em ordered triangle} if $I<J<K$.

We shall need notation for the edges of a given face in $\cH$.
For each $\beta\in\cH$, we let $\cE(\beta)\subset \cE $ be the subset
corresponding to the edges of $\beta$.  For example, if $\beta$ is the
triangular face $\{12,13,23\}$ of the octahedron, then
$\cE(\beta)=\{\{12,13\},\{12,23\},\{13,23\}\}$.  

%%%%%%%%%%%%%%%%%%%%
%%%
%%% The global resolution  
%%%
%%%%%%%%%%%%%%%%%%%
\section{The resolution $\tilde{X}$}\label{s:global}

\subsection{}\label{ss:tautological-bundles}

As the first step towards defining $\tilde{X}$, we establish a correspondence
between edges of the hypersimplices and certain $\proj^1$-bundles over
$X$.  For each nonempty subset $I\subset\mainset$, let
$F_I\rightarrow X$ be the pullback of the tautological $\card{I}$-plane
bundle on $\Grass_{\card{I}}$ via the composition
\[X\longrightarrow Y\stackrel{\pi_I}{\longrightarrow}\Grass_{\card{I}}.\]
Thus the fiber of $F_I$ over a point $p\in X$ can be identified with the
$k$-dimensional subspace $\pi_{I}(p)\subset\C^4$.  The incidence
conditions on $X$ imply that if $I\subset J$, then $F_I\subset F_{J}$ is a
subbundle.

For each $\alpha\in \cE$, let $F_{\alpha}$ be the
quotient $F_{I\cup J}/F_{I\cap J}$ where $\alpha=\{I,J\}$; this is a
rank $2$ vector bundle since $\card{I\cup J \smallsetminus I\cap J} =
2$. Let $P_{\alpha}$ be the projectivized bundle 
\[P_{\alpha}=\proj(F_{\alpha}).\]

\subsection{}
The bundle $P_{\alpha}$ has canonical sections
$u_{\alpha}^-,u_{\alpha}^+\colon X\rightarrow 
P_{\beta}$, defined geometrically as follows.  The fiber of
$P_{\alpha}$ over $p$ can be identified with the set of lines 
in the $2$-dimensional vector space $\pi_{I\cup J}(p)/\pi_{I\cap
J}(p)$.  We assume that $I<J$, and define
$u^-_{\alpha}(p)$ to be the line $\pi_{I}(p)/\pi_{I\cap J}(p)$ and
$u^+_{\alpha}(p)$ to be the line $\pi_{J}(p)/\pi_{I\cap J}(p)$.   
Since we will want to keep track of both sections
simultaneously, we introduce the product bundle 
\[E_{\alpha}=P_{\alpha}\times_{X}P_{\alpha},\]
and let $u_{\alpha}\colon X\rightarrow E_{\alpha}$ be the product
$u_{\alpha}^-\times u_{\alpha}^+$.  

\subsection{}
For each $\beta\in\cH$, let $E_{\beta}$ be the product bundle
\[E_{\beta}=\prod_{\alpha\in\cE(\beta)}E_{\alpha}.\]
This is a $(\proj^1\times\proj^1)^3$-bundle over $X$ when $\beta$ is a
triangular face; for the maximal faces, $P_{\Delta_k}$ is a
$(\proj^1\times\proj^1)^6$-bundle for $k=1,3$ and a
$(\proj^1\times\proj^1)^{12}$-bundle for $k=2$.  Let $p_{\beta}\colon
E_{\beta}\rightarrow X$ be the projection, and let $u_{\beta}\colon
X \rightarrow E_{\beta}$ be the section obtained by taking the product
of the sections $u_{\alpha}$ for all $\alpha\in\cE(\beta)$.

\subsection{}
We define the {\em ambient variety $E$} to be the product
bundle
\[E=\prod_{\beta\in\cH}E_{\beta}.\]
By Figure~\ref{labgeneric.fig}, there are $3$ maximal elements
of $\cH$ with $6$, $12$, and $6$ edges each, and there are $16$
triangular faces with $3$ edges each; thus, $E\rightarrow X$ is a
locally trivial bundle with fiber isomorphic to  
\[(\proj^1\times\proj^1)^6\times(\proj^1\times\proj^1)^{12}\times
(\proj^1\times\proj^1)^6\times
\left((\proj^1\times\proj^1)^3\right)^{16}.\] Let $p\colon
E\rightarrow X$ be the projection, and let $u\colon X\rightarrow E$ be
the section obtained by taking the product of the sections
$u_{\beta}$, $\beta\in\cH$.

\subsection{}\label{ss:beta-bu}

To build $\tilde{X}$, we keep track of ``limiting
configurations'' of the subspaces $\{\pi_I(p)\}$
as certain collections of them coincide.  The relevant collections
turn out to correspond to the faces $\cH$ of the hypersimplices.

For each $\alpha\in\cE$, let $D_{\alpha}\subset E_{\alpha}$ be 
the diagonal subbundle, and for each $\beta\in\cH$, let
$D_{\beta}\subset E_{\beta}$ be the subbundle  
\[D_{\beta}=\prod_{\alpha\in\cE(\beta)}D_{\alpha}.\]
The geometric significance of $D_{\beta}$ is that the set of points
$p\in X$ such that $u_{\beta}(p)\in D_{\beta}$ is precisely the set of $p$
such that $\pi_I(p)=\pi_J(p)$ for all $I,J\in\beta$. 

Let 
\[b_{\beta}\colon (E_{\beta})_{\#}\longrightarrow E_{\beta}\]
be the blowup of $E_{\beta}$ along $D_{\beta}$.
Since $E_{\beta}$ is locally trivial over $X$, as is the subbundle
$D_{\beta}$, the blowup $(E_{\beta})_{\#}$ is also locally trivial
over $X$; the fiber of this last bundle is isomorphic to
the blowup of $(\proj^1\times\proj^1)^n$ along
the product of diagonals (where $n=3$, $6$, or $12$ depending on $\beta$).  

\subsection{}\label{ss:compamb}

Since for any $p\in X^{\circ}$ the image $u_{\beta}(p)$ avoids the
blowup center $D_{\beta}$, we have a regular map  
\[b_{\beta}^{-1}\circ
u_{\beta}\colon X^{\circ}\longrightarrow (E_{\beta})_{\#}.\] 
We define
the {\em complete ambient variety $E_{\#}$} to be the product
\[E_{\#}=\prod_{\beta\in\cH}(E_{\beta})_{\#},\]
and let $b\colon E_{\#}\rightarrow E$ be the product of the blowup maps
$b_{\beta}$, $\beta\in\cH$.

\begin{definition}
Let $\tilde{X}^{\circ}$ be the image of the embedding 
\[X^{\circ}\longrightarrow E_{\#}\]
obtained by taking a product of the maps $b_{\beta}^{-1}\circ
u_{\beta}$ for all $\beta\in\cH$.  The {\em complete space of
tetrahedra}, denoted $\tilde{X}$, is the closure of
$\tilde{X}^{\circ}$ in $E_{\#}$. 
\end{definition}

The composition $p\circ b\colon E_{\#}\rightarrow X$ restricts to a
surjective birational morphism $\rho\colon \tilde{X}\rightarrow X$.

\begin{remark}
The map $b\colon E_{\#}\rightarrow E$ can be realized as an iterated blowup
along regularly embedded subschemes.  In this setting, the complete
space of tetrahedra $\tilde{X}$ is the (iterated) proper transform of
$u(X)\subset E$.  
\end{remark}

\subsection{}
Since the bundles $E_{\beta}$ are constructed from tautological
bundles, they admit natural $G$-actions lifting the action on $X$.
Since the diagonals are preserved by these actions, the blown-up
bundles $(E_{\beta})_{\#}$ also admit natural $G$-actions that lift
the action on $X$, and the blowdown maps $b_{\beta}$ are equivariant.
It follows that there are natural $G$-actions on $E$ and $E_{\#}$, and
that $b\colon E_{\#}\rightarrow E$ is equivariant.  One can check that
the section $u$ is also equivariant, and thus $\tilde{X}^{\circ}$ is
$G$-stable.  It follows that the action on $E_{\#}$ restricts to an
action on $\tilde{X}$ and that $\rho\colon \tilde{X}\rightarrow X$ is
$G$-equivariant.

Similar remarks apply to the $S_{4}$-action.  This action also lifts
to actions on $E$ and $E_{\#}$ that permute the various factors of
these product bundles.  The map $b\colon E_{\#}\rightarrow E$ and the
section $u\colon X\rightarrow E$ are both equivariant, so
$\tilde{X}^{\circ}$ is $S_4$-stable.  Hence, the $S_4$-action on
$E_{\#}$ restricts to an action on $\tilde{X}$, and $\rho\colon
\tilde{X}\rightarrow X$ is $S_4$-equivariant.

%%%%%%%%%%%
%
%  section four
%
%%%%%%%%%%%%

\section{The local variety $U$}\label{s:infinity}

\subsection{}\label{ss:schubertcell}
Let $\Flag$ be the flag variety of full flags
in $\C^4$, and let $V_*\in \Flag $ correspond to a flag
\[\{0 \} = V_{0} \subsetneq V_1\subsetneq V_2\subsetneq V_3 \subsetneq V^{4} = \C ^{4},\]
where $V_{k}$ is a subspace of dimension $k$.
Let $U(V_*)$ be the set of all $p\in X$ in general position to $V_*$.
That is, $U(V_*)$ consists of all $p$ such that the $\card{I}$-plane
$\pi_I(p)$ is transverse to $V_k$ for all $1\leq k\leq 3$ and all
proper nonempty subsets $I\subset\mainset$.   

The subset $U(V_*)$ can be described in terms of Schubert cells in the
factors $\Grass _{I}$ of $Y$ as follows.  For each $k=1,2,3$, let $U_k$
be the open cell in $\Grass_k$ consisting of $k$-planes in general
position to the fixed flag $V_*$.  For each proper
nonempty subset $I\subset\mainset$, let $U_I:=U_{\card{I}}$.
Then $\prod_I U_I$ is an open subvariety of $Y$ isomorphic to an
affine space.  The variety $U(V_*)$ is the intersection of this open
set and the subvariety $X$ of $Y$.  In particular, $U(V_*)$ is an
affine open subset of $X$.

\subsection{}
Let $E_{\#}|_{U(V_*)}$ be the restriction of the ambient bundle to $U
(V_{*})$.  Let $U^{\circ}(V_*)=U(V_*)\cap X^{\circ}$, and let
$\tilde{U}^{\circ}(V_*)$ be the image of $U^{\circ}(V_*)$ under the
embedding $U^{\circ}(V_*)\rightarrow E_{\#}|_{U(V_*)}$ of
\ref{ss:compamb}.  Take $\tilde{U}(V_*)$ to be the closure of
$\tilde{U}^{\circ}(V_*)$ in $E_{\#}|_{U(V_*)}$.

\begin{lemma}\label{lem:globalcharts}
The collection $\{U(V_*)\mid V_*\in\Flag\}$ (respectively,
$\{\tilde{U}(V_*)\mid V_*\in\Flag\}$) is an affine open cover of $X$
(resp., $\tilde{X}$).  The group $G$ acts transitively on both of these
covers.
\end{lemma}

To prove that $\tilde{X}$ is nonsingular, it suffices by
Lemma~\ref{lem:globalcharts} to prove that $\tilde{U}(V_*)$ is
nonsingular for one particular choice of the flag $V_{*}$.  We fix $V_*$ to
be the {\em flag at infinity}: 
\[E_4\subset E_{34}\subset
E_{234},\] 
and define $U$, $U^{\circ}$, and $\tilde{U}$ to be the varieties 
$U(V_*)$, $\tilde{U}^{\circ}(V_*)$, and $\tilde{U}(V_*)$
(respectively).

\subsection{}\label{ss:coordinates}
We put coordinates on $U$ using Pl\"ucker coordinates on the
Grassmannians $\Grass _{I}$.  For each $k=1,2,3$, we have the Pl\"ucker
embedding $\Grass_k\rightarrow\proj(\bigwedge^k\C^4)$ with
its usual coordinates $\{f_{I}\mid I\subset \mainset, \card{I}=k\}$.
The ratios $\{f_{I}/f_{I_{0}}\mid \card{I}=k \}$ provide
coordinates on $U_{k}$.  For any pair $I,J\subset \mainset $ with  
$\card{I}=\card{J}$, let $f_{I,J}$ be 
the regular function on $U$ defined by 
\[
f_{I,J} := \pi _{I}^{*} (f_{J}/f_{J_{0}}).
\] 
It is clear that these functions generate the ring
$\sO_{X}(U)$, and that $f_{I,J_{0}} = 1$.

\subsection{}
There is a more symmetric set of generators for $\sO_X(U)$, which
arises from the observation that $p\in U$ can be constructed from
functions on $\Flag$ and functions that measure the
``difference'' between the planes $\pi_I(p)$ and $\pi_J(p)$ for each edge
$\{I,J\}$ of the appropriate hypersimplex.

The functions on the flag variety are defined as follows.  There is a
natural map $X\rightarrow\Flag$ given by $p\mapsto \{\pi _{I_{0}}(p) \}$.
Let $U_{op}$ be the open cell in $\Flag$ consisting of flags in
general position to the flag at infinity.  This cell has local coordinates
\[f_2/f_1,\; f_3/f_1,\; f_4/f_1,\; f_{13}/f_{12},\; f_{14}/f_{12},\;
f_{124}/f_{123}.\] The corresponding functions
\[f_{1,2},\; f_{1,3},\; f_{1,4},\; f_{12,13},\; f_{12,14},\;
f_{123,124}\] on $U$ will be called {\em flag coordinates on $U$}.

\subsection{}\label{ss:sections} 
The functions on $U$ corresponding to edges in the
hypersimplices are easiest to describe using certain local
sections of the bundles of \ref{ss:tautological-bundles}.  

For each nonempty $I\subset\mainset$, let $\sF_I$ be the sheaf of
sections of the bundle $F_I\rightarrow X$.  We define local sections
$s_I\in\sF_I(U)$ as follows.  Let $k=\card{I}$.  Then for each $p\in
U$, the fiber of $F_I$ over $p$ can be identified with the
$k$-dimensional subspace $\pi_I(p)\subset\C^4$.  This subspace
intersects the subspace $V_{5-k}$ of our flag at
infinity in a 
$1$-dimensional subspace, and intersects $V_{4-k}$ in the
zero subspace.  It follows that there is a unique vector
$s_I(p)\in\pi_I(p)\cap V_{5-k}$ whose $k$th coordinate (with respect
to the standard basis) is $1$.  This defines the section $s_I\colon
U\rightarrow F_I|_{U}$.

In terms of Pl\"{u}cker coordinates, these sections can be expressed as
\[\begin{array}{lcrcrcrcrl}
s_i& = & e_1&+ & f_{i,2}\;e_2&+ & f_{i,3}\;e_3&+
 & f_{i,4}\;e_4&\\
s_{ij}&=& && e_2&+ &f_{ij,13}\;e_3&+
 &f_{ij,14}\;e_4&\\
s_{ijk}&=& & && &e_3&+  & f_{ijk,124}\;e_4&\\
s_{1234}&=& & && &&  & e_4 & .\end{array}
\]
A priori, these are all sections of the trivial bundle $U\times\C^4$,
but a simple verification shows that their images are contained in $F_I|_{U}$.  The following
lemma describes a crucial relation among these sections.  We omit the
straightforward proof.

\begin{lemma}\label{lem:ec-relations}
Let $k=1,2,3$.  For each edge $\alpha=\{I,J\}\in \cE_k$, we have
\[s_J-s_I=(f_{J,K}-f_{I,K})s_{I\cup J},
\]
where $K$ is the subset $2$, $13$, or $124$ depending on whether $k$
is $1$, $2$, or $3$, respectively.
\end{lemma}

For $k=1,2,3$ and each edge $\alpha=\{I,J\}\in \cE_k$ with $I<J$,
we define the {\em edge coordinate} $x_{\alpha}$ by
\[x_{\alpha}= f_{J,K}-f_{I,K},\]
where $K$ is the subset $2$, $13$, or $124$ depending on whether $k$
is $1$, $2$, or $3$, respectively.  

\begin{lemma}\label{lem:genfunc}
The ring $\sO_{X}(U)$ is generated by the flag coordinates and the
edge coordinates.
\end{lemma}

\begin{proof}{}
We show that the functions $f_{I,J}$ can be expressed in terms of the
flag coordinates and the $x_{\alpha}$'s.  The proof is by induction on $I$,
using the total ordering 
\[1\prec 12\prec 123\prec 2\prec 13\prec 3\prec 23\prec 124\prec 14\prec 4\prec 24\prec 134\prec 34\prec 234.\] 
The key property of the ordering $\prec$ is that for each $J\succ 123$, there
exists an edge $\{I,J\}\in\cE$ such that $I$, $I\cap J$, $I\cup J
\prec J$.
 
Using the formulas of \ref{ss:sections}, we can express the sections
$s_1$, $s_{12}$, and $s_{123}$ entirely in terms of the flag coordinates
(and the basis $e_1,e_2,e_3,e_4$).  Since $s_1$ determines the line
$\pi_1$, $s_1\wedge s_{12}$ determines the plane $\pi_{12}$, and
$s_1\wedge s_{12}\wedge s_{123}$ determines the $3$-plane $\pi_{123}$,
we can express all of the corresponding functions $f_{1,J}$,
$f_{12,J}$, and $f_{123,J}$ in terms of the flag coordinates.   

For the case $I=2$, since $s_2=s_1+x_{1,2}s_{12}$ (by
Lemma~\ref{lem:ec-relations}) and $s_2$ determines
$\pi_2$, we can express the functions $f_{2,J}$ in terms of
$x_{1,2}$ and the flag coordinates.  For $I=13$,
$s_{13}=s_{12}+x_{12,13}s_{123}$ and $s_1\wedge s_{13}$ determines
$\pi_{13}$, so we can express the functions $f_{13,J}$ in terms of 
$x_{12,13}$ and the flag coordinates.  Expressions for the remaining
functions are obtained similarly.
\end{proof}

\subsection{}
Define polynomial rings 
\begin{align*}
R_{op}&:= \C[f_{1,2},f_{1,3},f_{1,4},f_{12,13},f_{12,14},f_{123,124}],\\
R_{\cE} &:= \C[x_{\alpha}\mid \alpha\in\cE], 
\end{align*}
and let $\bbA_{op}=\Spec R_{op}$, $\bbA_{\cE}=\Spec R_{\cE}$.
Lemma~\ref{lem:genfunc} says that the natural homomorphism
$R_{op}\otimes R_{\cE}\rightarrow\sO_X(U)$ is surjective, so $U\subset
\bbA_{op}\times\bbA_{\cE}$.  We now describe the ideal that
set-theoretically cuts out $U$.  This ideal is generated by linear,
quadric, cubic, and quartic polynomials in the flag and edge
coordinates.

First consider the flag coordinates.  The map $X\rightarrow\Flag$ is
actually a locally trivial fibration, and our coordinates define a
trivialization over $U_{op}$.  Since $U_{op}$ is nonsingular and the
flag coordinates on $U$ are pulled back from a system of local
parameters on $U_{op}$, there are no relations among the flag
coordinates holding on $U$.

Now consider the edge coordinates.  Recall that a triple of edges
$(JK,IK,IJ)$ is an ordered triangle if $\{I,J,K\}$ is a triangular
face and $I<J<K$.

\begin{figure}[ht]
\begin{center}
\includegraphics[scale = .6]{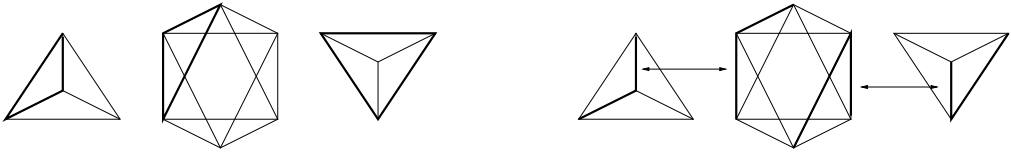}
\end{center}
\caption{\label{linquad.fig} Edges in the linear and quadric relations}
\end{figure}

\begin{figure}[ht]
\begin{center}
\includegraphics[scale = .6]{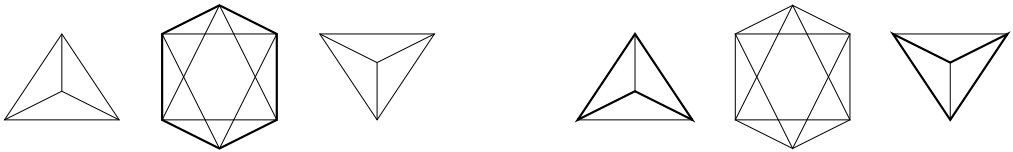}
\end{center}
\caption{\label{cubquar.fig} Edges in the cubic and quartic relations}
\end{figure}

\begin{lemma}\label{lem:all-equations}
The subvariety
$U\subset\bbA_{op}\times\bbA_{\cE}$ is defined set-theoretically by
the following polynomials:   
\begin{enumerate}
\item The linear functions
\[x_{\alpha_1}-x_{\alpha_2}+x_{\alpha_3},\]  
for all ordered triangles $(\alpha_1,\alpha_2,\alpha_3)$
(Figure~\ref{linquad.fig}).
\item The quadric functions
\[x_{\alpha_1}x_{\alpha_2^*}-x_{\alpha_2}x_{\alpha_1^*},\]
where $\alpha_1=\{i,j\}, \alpha_2=\{j,k\}, \alpha_1^*=\{ik,jk\},
\alpha_2^*=\{ij,ik\}$ or $\alpha_1=\{il,jl\}, \alpha_2=\{jl,kl\},
\alpha_1^*=\{ikl,jkl\}, \alpha_2^*=\{ijl,ikl\}$ (Figure~\ref{linquad.fig}).
\item The cubic functions
\[x_{\alpha_1}x_{\alpha_2}x_{\alpha_3}-
x_{\alpha_1^*}x_{\alpha_2^*}x_{\alpha_3^*},\]
where $\alpha_1=\{ij,il\}, \alpha_2=\{ik,kl\}, \alpha_3=\{jk,jl\},
\alpha_1^*=\{jk,kl\}, \alpha_2^*=\{ij,jl\}, \alpha_3^*=\{ik,il\}$
(Figure~\ref{cubquar.fig}).
\item The quartic functions
\[x_{\alpha_1}x_{\alpha_3}x_{\alpha_2^*}x_{\alpha_4^*}- 
x_{\alpha_2}x_{\alpha_4}x_{\alpha_1^*}x_{\alpha_3^*},\]
where $\alpha_1=\{i,j\}, \alpha_2=\{j,k\}, \alpha_3=\{k,l\},
\alpha_4=\{l,i\}, \alpha_1^*=\{ikl,jkl\}, \alpha_2^*=\{ijl,ikl\}, 
\alpha_3^*=\{ijk,ijl\}, \alpha_4^*=\{jkl,ijk\}$ (Figure~\ref{cubquar.fig}).
\end{enumerate}
\end{lemma}

\begin{proof}{}
The vanishing of the linear polynomials follows from the definition of the
edge coordinates.  The quadric relations follow from this definition
and Lemma~\ref{lem:ec-relations}.  The cubic (resp., quartic)
relations can be obtained from the quadric relations by
eliminating coordinates corresponding to edges in $\cE_1$ and $\cE_3$ 
(resp., $\cE_2$).  It follows that $U$ is a subvariety of the variety
defined by the given polynomials.  

To see that $U$ coincides with this variety, it suffices to show that
they are both irreducible and have the same dimension.  Since $U$ is
the closure of the connected $12$-dimensional nonsingular variety
$U^{\circ}$, it is a $12$-dimensional irreducible variety.  Let
$U_{inc}$ be the variety defined by the vanishing of the linear and
quadric polynomials.  A simple computation using Macaulay2
\cite{macaulay2} shows that $U_{inc}$ has three irreducible
components, two of which correspond to the ideals 
\[
\langle x_{\alpha} \mid \alpha \in \cE _{2}\rangle\quad
\text{and}\quad \langle x_{\alpha
}\mid \alpha \in \cE_{1}\cup \cE_{3}\rangle.
\]
The remaining component $U_{1}$ is $12$-dimensional.  Since $U$ contains
points where all $x_{\alpha }$ are nonzero, we have $U=U_{1}$.  
\end{proof}

\begin{figure}[ht]
\begin{center}
\includegraphics[scale = .6]{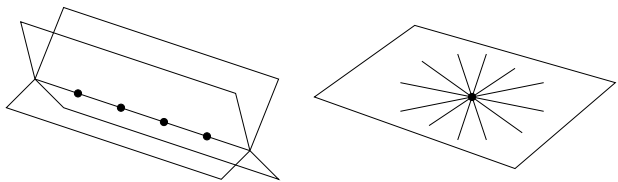}
\end{center}
\caption{\label{config.fig}Other components of the incidence variety.}
\end{figure}

\begin{remark}
Let $X_{inc}$ be the \emph{incidence variety} consisting of all $p\in
Y$ such that $\pi_I(p)\subset\pi_J(p)$ whenever $I\subset J$.  Then
one can show $U_{inc} = X_{inc}\cap \prod _{I}U_{I} $, where $\prod
_{I}U_{I}$ is the affine cell of \ref{ss:schubertcell}, although we
will not need this here.  Hence the linear and quadric relations
provide a very simple description of the incidence variety.  

In the case under study, Lemma \ref{lem:all-equations} shows that
$X_{inc}$ has two other components besides $X$, corresponding to the
following types of configurations in $\C ^{4}$:
\begin{enumerate}
\item Four $3$-planes containing a $2$-plane that contains four lines.
\item Six $2$-planes containing a line and contained in a $3$-plane.
\end{enumerate}
Figure~\ref{config.fig} shows general points in these
components as configurations in $\proj ^{3}$.  For the general case of
configurations in $\C ^{n}$, the
components of $X_{inc}$ are unknown.
\end{remark}

\section{The local resolution $\tilde{U}$}\label{s:infinity-tilde}

\subsection{}
In this section we describe an open subbundle $W$ (with affine space
fibers) of the restricted ambient bundle $E|_{U}$.  The section $u$
restricts to a section of $W$, thus we can apply the blowup
construction to $W$, obtaining an open subbundle $W_{\#}$ of
$E_{\#}|_U$ that contains the local resolution $\tilde{U}$.  To get
functions on $\tilde{U}$, we then show that $\tilde{U}$ is contained
in a certain closed subvariety of $W_{\#}$ defined in terms of
the edge coordinates.   

To describe $W$, we first need to describe some trivializations of the
various bundles restricted to $U$.

\subsection{}\label{ss:Falpha-triv}
For each $\alpha=\{I,J\}\in\cE$, let $\sF_{\alpha}$ be the sheaf of
sections of the quotient bundle $F_{\alpha}=F_{I\cup J}/F_{I\cap J}$.
Letting $\sF_{\alpha}^*$ be the dual module, and restricting to
$U$, we then have (by definition)  
\[F_{\alpha}|_{U}=\Spec_{\sO_X(U)} S_{\alpha}\]
where $S_{\alpha}$ is the symmetric $\sO_X(U)$-algebra
$\operatorname{Sym}\sF^*_{\alpha}(U)$.  

There are natural maps $\sF_I\rightarrow\sF_{\alpha}$,
$\sF_J\rightarrow\sF_{\alpha}$, and $\sF_{I\cup
J}\rightarrow\sF_{\alpha}$ induced by the corresponding inclusions of
$F_I$, $F_J$, and $F_{I\cup J}$ (respectively) into $F_{I\cup J}$.
With respect to these maps, the local sections $s_I\in\sF_I(U)$,
$s_J\in\sF_J(U)$, and $s_{I\cup J}\in\sF_{I\cup J}(U)$ all have images
in $\sF_{\alpha}(U)$, which we denote by $\overline{s}_I$,
$\overline{s}_J$, and $\overline{s}_{I\cup J}$, respectively.  It
follows from the explicit descriptions in \ref{ss:sections} that
$\overline{s}_I+\overline{s}_J$ and $\overline{s}_{I\cup J}$ are
nonzero and linearly independent.  

Since $F_{\alpha}$ is a rank-$2$
vector bundle, these sections determine dual sections $g_{\alpha},
h_{\alpha}\in \sF_{\alpha}^*(U)$, giving an isomorphism
%\[F_{\alpha}|_{U}\cong U\times\bbA^2,\]
\[
F_{\alpha}|_{U}\cong U\times\Spec \C [g_{\alpha },h_{\alpha }].
\]
Passing to the projectivized bundle
$P_{\alpha}=\proj(F_{\alpha})$ restricted to $U$, we obtain 
\[P_{\alpha}|_{U}=\Proj_{\sO_X(U)} S_{\alpha},\]
which becomes 
\[P_{\alpha}|_{U}\cong U\times\Proj \C [g_{\alpha }, h_{\alpha }].\]

\subsection{}
Let $D(g_{\alpha})\subset P_{\alpha}|_{U}$ be the divisor
determined by $g_{\alpha}$, and let $V_{\alpha}$ be the
corresponding open set $P_{\alpha}|_{U}-D(g_{\alpha})$.  Then $V_{\alpha}$
is an affine line bundle over $U$, and we have an isomorphism  
%\[V_{\alpha}\cong U\times\bbA^1\]
\[
V_{\alpha}\cong U\times \Spec \C [h_{\alpha }/g_{\alpha }].
\]

\subsection{}
For each $\alpha\in\cE$, let $W_{\alpha}$ be the product
$V_{\alpha}\times_U V_{\alpha}$.  To distinguish the two factors of
$W_{\alpha}$, we denote the first by $V_{\alpha}^-$ and the second by
$V_{\alpha}^+$.  For any $\beta\in\cH$, let
$W_{\beta}$ be the product of the $W_{\alpha}$ as $\alpha$ runs
over all edges in $\cE(\beta)$.  Finally, let $W$ be the product of
the $W_{\beta}$ as $\beta$ ranges over all faces in $\cH$.  Combining
the previous sections, we then have the following:

\begin{lemma}\label{lem:ambient-variables}
Let $R$ be the polynomial ring 
\[
R := \C [x_{\alpha ,\beta }^{+}, x_{\alpha, \beta }^{-} \mid
\beta\in\cH, \alpha\in\cE(\beta)],
\]
and let $\bbA_{amb}=\Spec R$.  Then $W$ is an open subbundle of
$E|_{U}$, and is isomorphic to the product $U\times\bbA_{amb}$. 
The indeterminate $x_{\alpha,\beta}^-$ (resp., $x_{\alpha,\beta}^+$)
corresponds to the function $h_{\alpha}/g_{\alpha}$ on the factor
$V_{\alpha}^-$ (resp., $V_{\alpha}^+$) of $W_{\beta}$.  
\end{lemma}  

\subsection{}\label{ss:proj-beta} 
We now restrict the blowup construction of \ref{ss:beta-bu} to the
affine space bundle $W\rightarrow U$.  For each $\beta\in\cH$, we let
$(W_{\beta})_{\#}$ be the blowup of $W_{\beta}$ along the product of
diagonals $D_{\beta}\cap W_{\beta}$.  In terms of the coordinates in
Lemma~\ref{lem:ambient-variables}, the ideal defining this product of
diagonals is $\langle x_{\alpha,\beta}^+-x_{\alpha,\beta}^- \mid
\alpha\in\cE(\beta)\rangle$.  Thus, if $R_{\beta}$ is the polynomial
ring $\C[y_{\alpha,\beta}\mid \alpha\in\cE(\beta)]$ and $\proj_{\beta}
:= \Proj R_{\beta}$, we have a closed embedding
\[(W_{\beta})_{\#}\longrightarrow W_{\beta}\times\proj_{\beta},\]
where the ideal defining the image is
\[\bigl\langle y_{\alpha,\beta}(x^+_{\alpha^*,\beta}-x^-_{\alpha^*,\beta})-
y_{\alpha^*,\beta}(x^+_{\alpha,\beta}-x^-_{\alpha,\beta})\bigm | \alpha,\alpha^*\in\cE(\beta)\bigr\rangle.
\]

\subsection{}\label{ss:bu-embedding}
Now we take products over $\cH $.
Let $W_{\#}$ be the product over $U$ of the blowups
$\{(W_{\beta})_{\#} \mid \beta\in\cH\}$.  We then have an embedding
\[W_{\#}\longrightarrow W\times\prod_{\beta\in\cH}\proj_{\beta},\] 
and the blowdown map $b\colon W_{\#}\rightarrow W$ is simply the restriction
of the projection to the first factor.  Moreover, in terms of 
coordinates in Lemma~\ref{lem:ambient-variables}, the image of this
embedding is cut out by the multihomogeneous polynomials 
\[y_{\alpha,\beta}(x^+_{\alpha^*,\beta}-x^-_{\alpha^*,\beta})-
y_{\alpha^*,\beta}(x^+_{\alpha,\beta}-x^-_{\alpha,\beta})\]
for all $\beta\in\cH$ and $\alpha,\alpha^*\in\cE(\beta)$. 

\subsection{}
We now consider the subvariety $\tilde{U}\subset W_{\#}$.  This is, by 
definition, the closure of $\tilde{U}^{\circ}=b^{-1}(u(U^{\circ}))$. 

\begin{lemma}\label{lem:U-in-W}
The image of the section $u\colon U\rightarrow E|_{U}$ is contained in
the open subvariety $W$.  In terms of the coordinates of
Lemma~\ref{lem:ambient-variables} the section $u\colon U\rightarrow W$
is defined by the $\sO_X(U)$-module homomorphism defined by
$x_{\alpha,\beta}^{\pm}\mapsto\pm x_{\alpha}$.    
\end{lemma}

\begin{proof}{}
Let $\alpha=\{I,J\}\in\cE$ with $I<J$.  Then the section
$u_{\alpha}^-$, (respectively, $u_{\alpha}^+$) is defined, at
each point $p\in U$, to be the linear span of the nonzero vector
$\overline{s}_I(p)$ (resp., $\overline{s}_J(p)$).  With $g_{\alpha}$
as in \ref{ss:Falpha-triv}, it follows from the formulas in  
\ref{ss:sections} that $g_{\alpha}(\overline{s}_I)$ and 
$g_{\alpha}(\overline{s}_J)$ are nonzero on $U$.  Hence the image 
$u_{\alpha}(U)$ is contained in $W_{\alpha}$.  Taking suitable
products of these sections, we then have $u(U)\subset W$.    

Using the various bundle trivializations above, one can show that the
section $u$ is given by $x_{\alpha,\beta}^-\mapsto
h_{\alpha}(\overline{s}_I)/g_{\alpha}(\overline{s}_I)$ and
$x_{\alpha,\beta}^+\mapsto
h_{\alpha}(\overline{s}_J)/g_{\alpha}(\overline{s}_J)$).  By
Lemma~\ref{lem:ec-relations}, we have
\[\overline{s}_I-f_{I,K}\overline{s}_{I\cup J}=
\overline{s}_J-f_{J,K}\overline{s}_{I\cup J}.\] 
Applying $g_{\alpha}$ to this equation, and using
$g_{\alpha}(\overline{s}_I+\overline{s}_J)=1$, we have
$g_{\alpha}(\overline{s}_I)=g_{\alpha}(\overline{s}_J)=1/2$.  Applying
$h_{\alpha}$ to this equation, and using
$h_{\alpha}(\overline{s}_I+\overline{s}_J)=0$, we have
$h_{\alpha}(\overline{s}_J)=-h_{\alpha}(\overline{s}_I)=
(f_{J,K}-f_{I,K})/2=x_{\alpha}/2$.  Thus $u\colon U\rightarrow W$ is given
by $x_{\alpha,\beta}^{\pm}\mapsto\pm x_{\alpha}$.  
\end{proof} 

\subsection{}\label{ss:W-affine}
Combining the trivialization of Lemma~\ref{lem:ambient-variables} with
the embedding of Lemma~\ref{lem:all-equations}, we can view $W$ as a
subvariety of the affine space $\bbA_{op}\times\bbA_{\cE}\times
\bbA_{amb}$.  It follows from Lemma~\ref{lem:U-in-W} that the section
$u\colon U\rightarrow W$ is the restriction of the inclusion  
$\bbA_{op}\times\bbA_{\cE}\rightarrow
\bbA_{op}\times\bbA_{\cE}\times\bbA_{amb}$  
defined by $x_{\alpha,\beta}^{\pm}\mapsto\pm x_{\alpha}$; thus, $u(U)$
is defined set theoretically by the polynomials of
Lemma~\ref{lem:all-equations} together with the linear polynomials 
\[x_{\alpha,\beta}^+-x_{\alpha}\hspace{.5in}\mbox{and}\hspace{.5in} 
x_{\alpha,\beta}^-+x_{\alpha}\]
for all $\beta\in\cH$ and $\alpha\in\cE(\beta)$. 

\subsection{}\label{ss:main-embedding}
By combining the embedding of \ref{ss:bu-embedding} with
\ref{ss:W-affine}, the blowup $W_{\#}$ (and hence $\tilde{U}$) can be
regarded as a subvariety of  
\[\bbA_{op}\times\bbA_{\cE}\times
\bbA_{amb}\times\prod_{\beta\in\cH}\proj_{\beta}.\]   
It follows from the relations in \ref{ss:W-affine} that
$\tilde{U}^{\circ}$ (and hence its closure $\tilde{U}$) will be
contained in the subvariety defined by $x_{\alpha,\beta}^{\pm}=\pm
x_{\alpha}$ for all $\beta\in\cH$ and $\alpha\in\cE(\beta)$.  Since the
projection  
\[\bbA_{op}\times\bbA_{\cE}\times 
\bbA_{amb}\times\prod_{\beta\in\cH}\proj_{\beta}\longrightarrow
\bbA_{op}\times\bbA_{\cE}\times\prod_{\beta\in\cH}\proj_{\beta}\]
is an isomorphism when restricted to this subvariety, the further
restriction to $\tilde{U}$ defines a closed embedding
\[\tilde{U}\longrightarrow
\bbA_{op}\times\bbA_{\cE}\times\prod_{\beta\in\cH}\proj_{\beta}.\]
The image of this embedding is cut out by the
polynomials of Lemma~\ref{lem:all-equations} together with the 
multihomogeneous polynomials
\[y_{\alpha,\beta}x_{\alpha^*}-y_{\alpha^*,\beta}x_{\alpha}, \quad
\beta\in\cH, \quad \alpha,\alpha^*\in\cE(\beta).\]

\section{The core $\core$}\label{s:core}

\subsection{}
In this section we use the embedding of \ref{ss:main-embedding} to
show that $\tilde{U}$ is isomorphic to a $9$-dimensional vector bundle
over a certain $3$-dimensional multi-projective variety.  

\begin{definition}
Let $\tilde{U}\rightarrow\prod_{\beta}\proj_{\beta}$ be the
composition of the embedding of \ref{ss:main-embedding} with the
projection to the projective spaces.  The image of this map will be
called the {\em core}, and denoted $\core$.  We let
$\eta:\tilde{U}\rightarrow\core$ denote the induced map.
\end{definition}

\subsection{}
For each $k=1,2,3$, we consider the projection
$\core\rightarrow\proj_{\Delta_k}$, and let $L_k\rightarrow \core$ be the
pullback of the tautological line bundle.  Since
$\proj_{\Delta_k}=\Proj R_{\Delta_k}$ (see \ref{ss:proj-beta}), $L_k$ is
naturally a subvariety of $\Spec R_{\Delta_k}\times \core$.
By identifying the ring $R_{\cE}$ with $R_{\Delta_1}\otimes
R_{\Delta_2}\otimes R_{\Delta_3}$ (via $y_{\alpha,\Delta_k}\mapsto
x_{\alpha}$), we can identify the $3$-dimensional bundle
$L_1\times_{\core} L_2\times_{\core} L_3$ with a 
subvariety of $\bbA_{\cE}\times\core$.  Let $N\rightarrow \core$ be the
$9$-dimensional vector bundle obtained by taking the product (over $\core$)
of the trivial bundle $\bbA_{op}\times\core$ and the bundle
$L_1\times_{\core} L_2\times_{\core} L_3$.  There is a natural
embedding 
\[N\longrightarrow
\bbA_{op}\times\bbA_{\cE}\times\prod_{\beta\in\cH}\proj_{\beta},\]  
and it follows from the equations in \ref{ss:main-embedding} that
$\tilde{U}$ is contained in the image.  Thus, we have an embedding  
\[\tilde{U}\longrightarrow N\]
whose composition with the bundle projection to $\core$ coincides with
the map $\eta$. 

\subsection{}\label{ss:B-action}
To prove that $\tilde{U}$ coincides with $N$, we use the $G$-action on
$X$.  The stabilizer of the flag at infinity 
is the subgroup $B$ of $G$ consisting of lower triangular matrices.
The group $B$ acts on the varieties $U^{\circ}$ and $U$ by the usual
action on the Pl\"{u}cker coordinates.  In this section, we describe
a $B$-action on the bundle $N\subset\bbA_{op}\times\bbA_{\cE}\times\core$,
with the property that the embedding $\tilde{U}^{\circ}\rightarrow N$ is
$B$-equivariant.  

The action on $\core$ is trivial.  The action on $\bbA_{op}$ is the usual
action of the Borel on the corresponding big cell $U_{op}$ in the flag
variety.  The action on $\bbA_{\cE}$ is given in terms of the characters
$t_k:B\rightarrow\C^{\times}$ defined by $t_1(b)=b_{22}/b_{11}$,
$t_2(b)=b_{33}/b_{22}$, and $t_3(b)=b_{44}/b_{33}$ where $b$ is the
matrix $(b_{ij})$.  For each $\alpha\in\cE_k$, the action of $b$ on
$x_{\alpha}$ is then the diagonal action $x_{\alpha}\mapsto t_k(b)
x_{\alpha}$.  It is clear that $N$ is a $B$-stable subvariety of
$\bbA_{op}\times\bbA_{\cE}\times\core$. 

\begin{lemma}\label{lem:bundle-over-core}
The embedding $\tilde{U}\rightarrow N$ is a $B$-equivariant isomorphism.
\end{lemma}

\begin{proof}{}
Since the flag coordinates on $U$ are pulled back from the
coordinates on the flag variety, the composition
$\tilde{U}^{\circ}\rightarrow N\rightarrow\bbA_{op}$ is equivariant.
An explicit calculation using the sections of \ref{ss:sections} and
the definition of the edge coordinates shows
that the composition $\tilde{U}^{\circ}\rightarrow
N\rightarrow\bbA_{\cE}$ is equivariant.  And finally, since for each
$\beta\in\cH$, the group $B$ acts via the same character on
$x_{\alpha}$, for all $\alpha\in\cE(\beta)$, the induced action on each
$\proj_{\beta}$ will be trivial.  It follows that
$\tilde{U}^{\circ}$ embeds equivariantly into $N$; hence, so does its
closure.

To see that the embedding is an isomorphism, we let
${\core}^{\circ}=\eta(\tilde{U}^{\circ})$.  Since $\tilde{U}$ is the 
closure of $\tilde{U}^{\circ}$ in $W_{\#}$,  ${\core}^{\circ}$ is dense
in $\core$.  It follows from the description of the $B$-action that
the unipotent subgroup of $B$ acts freely and transitively on
$\bbA_{op}$, and that the diagonal subgroup of $B$ acts fiberwise on
the product of the complements of the zero sections in
$L_1\times_{\core} L_2\times_{\core} L_3$.  Thus $B$ acts with dense
orbit on  each fiber of $N$.  Since each $x_{\alpha}$ is nonzero on
the image of $\tilde{U}^{\circ}$ in $N$, the image of
$\tilde{U}^{\circ}$ intersects this $B$-orbit for every fiber of 
$N|_{{\core}^{\circ}}\rightarrow {\core}^{\circ}$.  It follows that
the image of $\tilde{U}^{\circ}$ is dense in $N$, so the image of its
closure $\tilde{U}$ coincides with $N$.  
\end{proof}

\subsection{}
By Lemma~\ref{lem:bundle-over-core}, we know that $\tilde{U}$ is
isomorphic to a vector bundle over $\core$, hence $\tilde{U}$ will be
smooth if and only if $\core$ is smooth.  We next show that
nonsingularity of the core $\core $ follows from nonsingularity along a
certain subvariety, called the locus of \emph{special points}.  We
begin with some notation. 

For any point $\tilde{p}\in \tilde{X}$, let $p$ be its image
in $X$.  We define the {\em number of $k$-planes in
$\tilde{p}$} by
\[
n_k(\tilde{p}):= \operatorname{Card}\bigl\{\pi_I(p)\bigm| I\subset\{1,2,3,4
\},\;\card{I}=k\bigr\}.
\]    
The point $\tilde{p}$ is called {\em split}
(resp. {\em minimally split}) if $n_k(\tilde{p})\geq 2$
(resp. $n_k(\tilde{p})=2$) for all $k\leq 3$.  A point $z\in \core $
is called \emph{special} if there exists a minimally split point
$\tilde{p}\in \tilde{U}$ with $\eta (\tilde{p}) = z$.  We let
$\core_{sp}\subset Z$ be the subvariety of special points.

\begin{proposition}\label{prop:toroidalfibers}
The split points are open in each fiber of
$\eta \colon \tilde{U}\rightarrow\core$.
\end{proposition}

\begin{proof}
For any fiber, we can choose $\tilde{p}\in \tilde{U}$ whose $B$-orbit
is open in that fiber.  The description of the $B$-action in
\ref{ss:B-action} therefore implies that for any $k\leq 3$, there will
be some $\alpha\in\cE_k$ such that $x_{\alpha}\not =0$ on the image of
$\tilde{p}$ in $N$.  But $x_{\alpha}\not =0$ implies
$\pi_I(p)\neq\pi_J(p)$, where $\alpha=\{I,J\}$.  Therefore $\tilde{p}$
is split, and the result follows since $n_{k}$ is constant on
$B$-orbits. 
\end{proof} 

\begin{proposition}\label{prop:2-flat-prop}
If $\tilde{p}\in \tilde{U}$ is split, then
$\closure{G\cdot\tilde{p}}\cap \tilde{U}$ contains a minimally split
point.
\end{proposition}

To prove the proposition we require some lemmas.  

\begin{lemma}\label{3flat-lemma}
Let $k\in\{1,2,3\}$ and let $F_{1}, F_{2}, F_{3}\in \Grass_k$ be three
distinct points.  Then there exists a one-parameter subgroup
$\mu:\C^{\times}\rightarrow G$ such that
\[\lim_{t\rightarrow 0}\mu(t)\cdot F_{2}=F_{1}\hspace{.4in}\mbox{and}\hspace{.4in}\lim_{t\rightarrow 0}\mu(t)\cdot F_{3}\neq F_{1}.\]
\end{lemma}

\begin{proof}
We can find a subspace $F_{4}\subset \C ^{4}$ such that $F_{4}\oplus
F_{1} = F_{4}\oplus F_{2} = \C ^{4}$, and such that $\dim F_{4}\cap
F_{3}>0$.  Then for $\mu $ we can take any one-parameter subgroup that
scales in $F_{1}$ with a negative weight and scales in $F_{4}$ with a
positive weight.
\end{proof}

\begin{lemma}\label{lem:t-lift1}
Let $\tilde{p},\tilde{s}\in \tilde{U}$ with $\tilde{p}$ split and
$\tilde{s}\in \closure{G\cdot \tilde{p}}$.  Then there exists a split
point $\tilde{r}\in \closure{G\cdot \tilde{p}}\cap \tilde{U}$ such
that $\tilde{s}\in \closure{B\cdot \tilde{r}}$.  Moreover, if $n_{k}
(\tilde{s}) > 1$, then $n_{k}(\tilde{s})=n_{k}(\tilde{r})$.
\end{lemma}

\begin{proof}
Let $W$ be the set of split points in $G\cdot \tilde{p}$.  This set is
open, and thus $\tilde{s}$ is in its closure.  But since the $B$-orbit
of any point in $W$ lies in $W$, the entire fiber $\eta ^{-1} (\eta
(\tilde{s}))$ must also be in the closure of $W$.  Letting $\tilde{r}$
be a split point in this fiber completes the proof of the first
statement.  For the second statement, a look at the $B$-action shows
that passing to a point in $\tilde{U}$ that is in an orbit closure
either preserves $n_{k}$ or drops it down to $1$.
\end{proof}

\begin{proof}[Proof of Proposition \ref{prop:2-flat-prop}]
We use the lemmas above to collapse the configuration associated to
$\tilde{p}$ so that only two subspaces of each dimension remain.  We
use implicitly that in passing to a point in the closure of a
$G$-orbit, the number of planes in any given dimension cannot
increase.  

We begin with the subspaces of dimension $1$.  By assumption,
$n_{1}(\tilde{p})>1$.  If $n_{1}(\tilde{p})>2$, then we can use Lemma
\ref{3flat-lemma} to find $\tilde{o}\in \closure{G\cdot\tilde{p}}$
such that $n_{1}(\tilde{p})>n_{1} (\tilde{o})\geq 2$.  The orbit
$G\cdot \tilde{o}$ must lie in $\closure{G\cdot\tilde{p}}$, and since
$G$ acts transitively on our charts that cover $\tilde{X}$ (see
\ref{lem:globalcharts}), we can find a $G$-translate $\tilde{s}$ of $\tilde{o}$
such that $\tilde{s}\in \closure{G\cdot\tilde{p}}\cap \tilde{U}$.
Since the functions $n_{k}$ are constant on $G$-orbits, we have
$n_{1}(\tilde{s}) = n_{1} (\tilde{o})$.  

Lemma~\ref{lem:t-lift1} implies that we can find a split point
$\tilde{r}\in \closure{G\cdot\tilde{p}}$ with $\tilde{s}\in
\closure{B\cdot\tilde{r}}\cap \tilde{U}$, and such that
$n_{1}(\tilde{r}) = n_{1}(\tilde{s})$ is $\geq 2$ and is $<n_{1}
(\tilde{p})$.  Since any point in $\closure{G\cdot\tilde{r}}$ is also
in $\closure{G\cdot\tilde{p}}$, we can repeat this procedure until we
find a split point $\tilde{p}_1\in \closure{G\cdot\tilde{p}}\cap
\tilde{U}$ with $n_{1} (\tilde{p}_{1}) = 2$.

Now we induct on $k$ to produce points $\tilde{p}_{2}$ and
$\tilde{p}_{3}$.  The key point is that we can apply the lemmas to
reduce $n_{k}$ while preserving $n_{l}>1$ for $l\not =k$.  Since the
final point $\tilde{p}_{3}\in \closure{G\cdot\tilde{p}}\cap \tilde{U}$
is minimally split, this completes the proof. 
\end{proof}

\subsection{}
Propositions~\ref{prop:toroidalfibers} and~\ref{prop:2-flat-prop}
imply that nonsingularity of $\tilde{U}$ follows from nonsingularity
at minimally split points.  Since $\tilde{U}$ is a vector bundle over
$\core $, nonsingularity at the minimally split points follows from
the nonsingularity of $\core _{sp}$.

%%%%%%%%%%%%%%%%%%%%
%
%   Nonsingularity 
%
%%%%%%%%%%%%%%%%%%%%
\section{Nonsingularity}\label{s:nonsing}

\subsection{}
Recall that the core $\core $ is a subvariety of 
\[
\prod _{\beta \in \cH } \proj _{\beta } \cong \proj ^{5}\times \proj
^{11}\times \proj ^{5}\times (\proj ^{2})^{4}\times (\proj ^{2}
)^{8}\times (\proj ^{2})^{4},
\]
where the index set $\cH $ corresponds to faces of the hypersimplices
of dimension $\geq 2$.  Each factor $\proj _{\beta }$ has homogeneous
coordinates $\{y_{\alpha ,\beta }\mid \alpha \in \cE (\beta ) \}$,
corresponding to the edges of the face $\beta $.  By combining the
polynomials of \ref{ss:main-embedding} with the polynomials defining
$U$ in \ref{lem:all-equations}, we obtain polynomials defining $\core$.  

\begin{lemma}\label{lem:core-equations}
The subvariety $Z$ of $\prod_{\beta}\proj_{\beta}$ is defined
set-theoretically by the following multihomogeneous polynomials:
\begin{enumerate}
\item The linear polynomials
\[y_{\alpha_1,\beta}-y_{\alpha_2,\beta}+y_{\alpha_3,\beta},\]
where $\alpha_1,\alpha_2,\alpha_3$ are as in
Lemma~\ref{lem:all-equations}(item 1) and $\beta\in\cH$ is such that
$\alpha_1,\alpha_2,\alpha_3\in\cE(\beta)$.
\item The quadric polynomials
\[y_{\alpha_1,\beta}y_{\alpha_2,\beta^*}-
y_{\alpha_2,\beta}y_{\alpha_1,\beta^*},\] 
where $\alpha_1$ and $\alpha_2$ are any two edges that share a vertex
and $\beta,\beta^*\in\cH$ are such that $\alpha_1,\alpha_2\in\cE(\beta)$ and
$\alpha_1,\alpha_2\in\cE(\beta^*)$. 
\item The quadric polynomials
\[y_{\alpha_1,\beta}y_{\alpha_2^*,\beta^*}-
y_{\alpha_2,\beta}y_{\alpha_1^*,\beta^*},\] 
where $\alpha_1,\alpha_2,\alpha_1^*,\alpha_2^*$ are as in
Lemma~\ref{lem:all-equations}(item 2) and $\beta,\beta^*\in\cH$ are
such that $\alpha_1,\alpha_2\in\cE(\beta)$
and $\alpha_1^*,\alpha_2^*\in\cE(\beta^*)$.
\item The cubic polynomials
\[y_{\alpha_1,\beta}y_{\alpha_2,\beta}y_{\alpha_3,\beta}-
y_{\alpha_1^*,\beta}y_{\alpha_2^*,\beta}y_{\alpha_3^*,\beta}\]
where $\alpha_1,\alpha_2,\alpha_3,\alpha_1^*,\alpha_2^*,\alpha_3^*$ are
as in Lemma~\ref{lem:all-equations}(item 3) and $\beta$ is the hypersimplex
$\Delta_2$. 
\item The quartic polynomials
\[y_{\alpha_1,\beta}y_{\alpha_3,\beta}y_{\alpha_2^*,\beta^*}
y_{\alpha_4^*,\beta^*}-  
y_{\alpha_2,\beta}y_{\alpha_4,\beta}y_{\alpha_1^*,\beta^*}
y_{\alpha_3^*,\beta^*}\]
where $\alpha_1,\alpha_2,\alpha_3,\alpha_4,\alpha_1^*,\alpha_2^*,
\alpha_3^*,\alpha_4^*$ are as in Lemma~\ref{lem:all-equations}(item 4)
and $\beta,\beta^*$ are the hypersimplices $\Delta_1,\Delta_3$,
respectively. 
\end{enumerate}
\end{lemma}

\subsection{}
\begin{figure}[ht]
\begin{center}
\includegraphics[scale = .4]{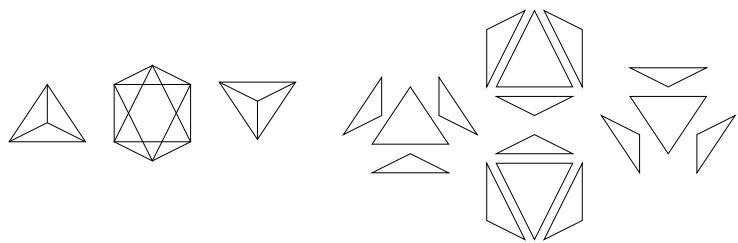}
\end{center}
\caption{\label{core.fig} The graph $\Gamma$}
\end{figure}

We represent points in $\prod _{\beta }\proj _\beta $ combinatorially
using the graph $\Gamma $ in Figure~\ref{core.fig}.  The edges of
$\Gamma $ are in bijection with the variables $y_{\alpha ,\beta }$,
and we encode a point in $\prod _{\beta }\proj _\beta$ by assigning
values to the edges modulo 
$(\C^{\times })^{19}$ (since $\Gamma $ has $19$ connected components).
It will be convenient to abuse language slightly by identifying points
in $\prod _{\beta }\proj _{\beta }$ with $\Gamma $.  In doing so we
shall always assume that values have been assigned to the variables
$y_{\alpha ,\beta }$.

Let $T$ and $T^*$ be two triangular subgraphs of $\Gamma $.  Then $T$
(respectively, $T^*$) corresponds to a choice of an ordered triangle
and a choice of a face in $\cH$ to designate the component of $\Gamma$ that 
contains $T$ (resp., $T^*$).  We shall say that $T$ and $T^*$
are {\em related} if one of the following two conditions holds:
\begin{enumerate}
\item The two ordered triangles for $T$ and $T^*$ coincide, or 
\item The faces in Figure~\ref{labgeneric.fig} 
that correspond to the two ordered triangles for $T$ and $T^*$ 
are in adjacent hypersimplices and one is a $180^{\circ}$ rotated copy
of the other. 
\end{enumerate}
In either case, there is a natural correspondence between the three
edges of $T$ and the three edges of $T^*$, and we say that $T$ and
$T^*$ have the \emph{same shape} if the corresponding triples of
values are proportional.  Our calculations will involve the use
of this notion together with setting various $y_{\alpha ,\beta }$ to
$0$; we indicate the latter by marking in bold the corresponding edge
of $\Gamma $.  

As a first step towards showing that special points are nonsingular,
we consider the equations in Lemma~\ref{lem:core-equations} and their
meaning in the context of Figure~\ref{core.fig}.   

\begin{figure}[ht]
\begin{center}
\psfrag{0}{$\scriptscriptstyle{0}$}
\includegraphics[scale = .3]{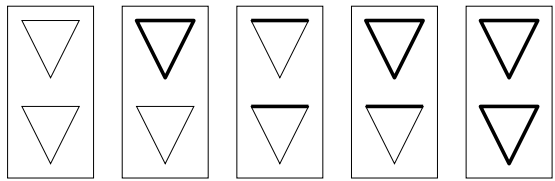}
\end{center}
\caption{\label{vert-rels.fig}Possible related triangles.}
\end{figure}

\begin{lemma}  
Suppose that $\Gamma $ represents a point in $Z$.  Then any
two related triangular subgraphs $T$ and $T^*$ have the same shape and
they must appear as one of the five possibilities shown in
Figure~\ref{vert-rels.fig}. 
\end{lemma}

\begin{proof}
The quadric relations in Lemma~\ref{lem:core-equations} imply that
related triangles will have the same shape.  Using this fact and the
linear relations of Lemma~\ref{lem:core-equations} it is then easy to
verify that the only combinations of zero values for such a pair are
those shown in Figure~\ref{vert-rels.fig}. 
\end{proof}

\begin{figure}[ht]
\begin{center}
\includegraphics[scale = .4]{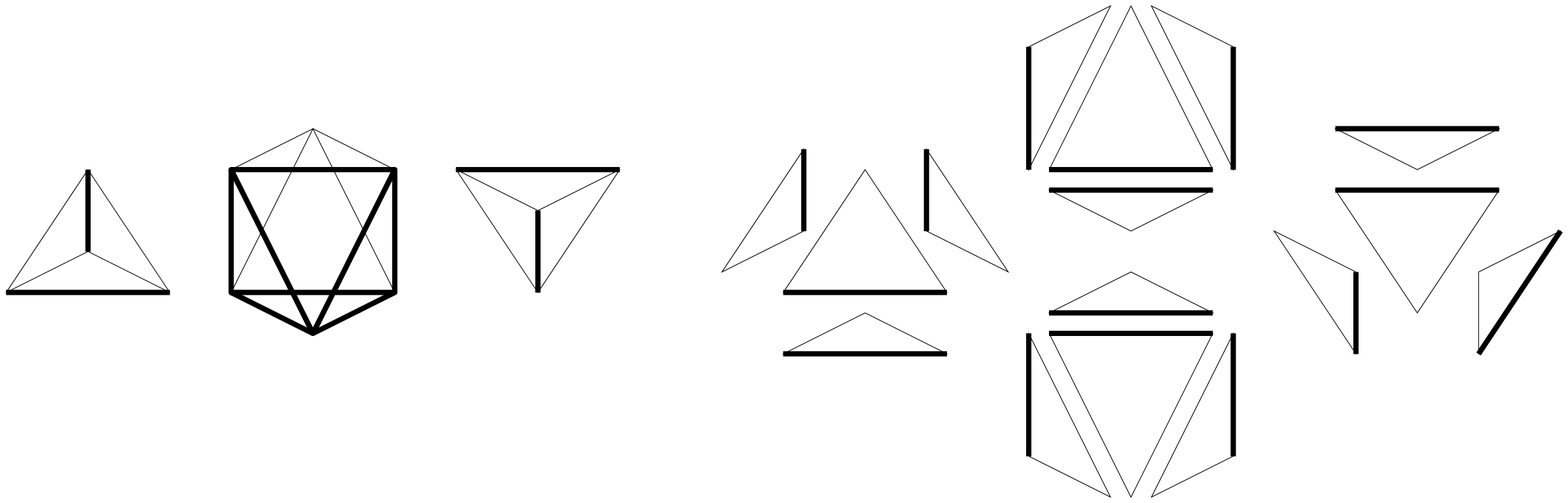}
\end{center}
\caption{\label{coreDDE.fig}$DDE\quad (6)$.}
\end{figure}

\begin{figure}[ht]
\begin{center}
\includegraphics[scale = .4]{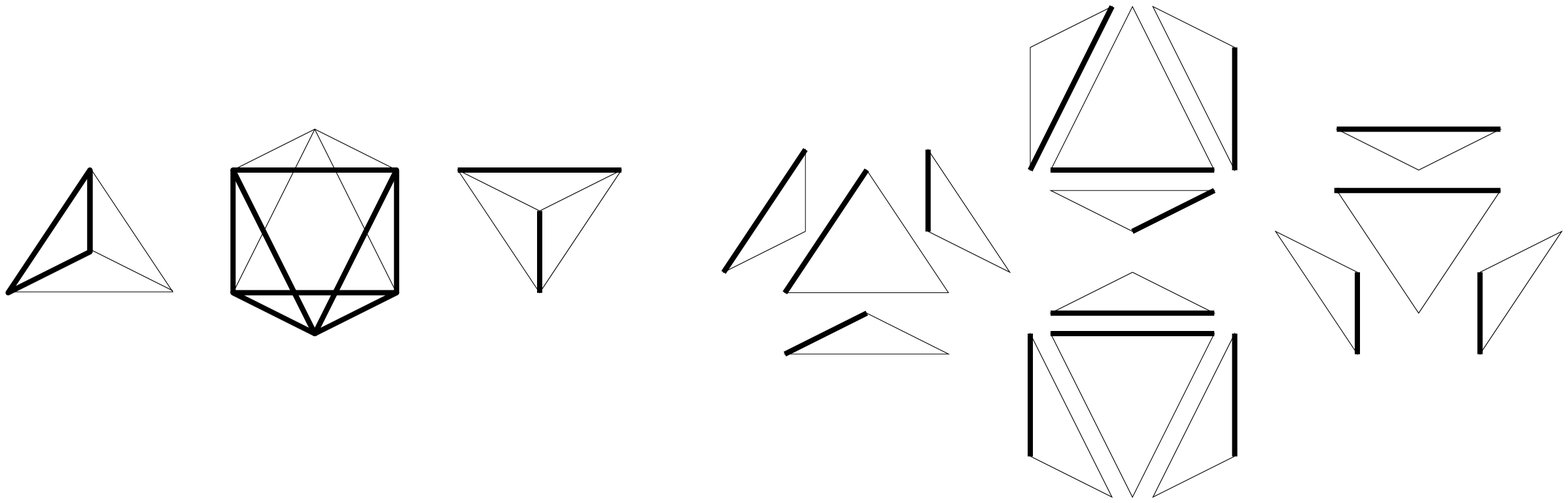}
\end{center}
\caption{\label{coreCDE.fig}$CDE\quad (12+12)$.}
\end{figure}

\begin{figure}[ht]
\begin{center}
\includegraphics[scale = .4]{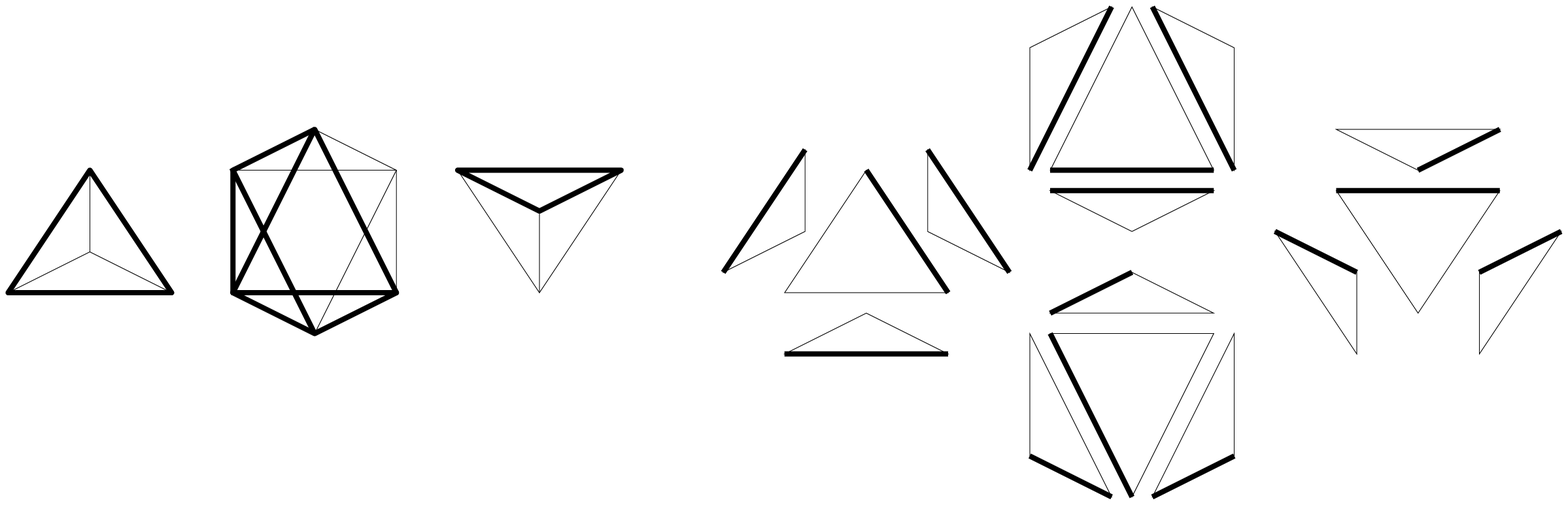}
\end{center}
\caption{\label{coreCCE.fig}$CC^{*}E \quad (24)$.}
\end{figure}

\begin{figure}[ht]
\begin{center}
\includegraphics[scale = .4]{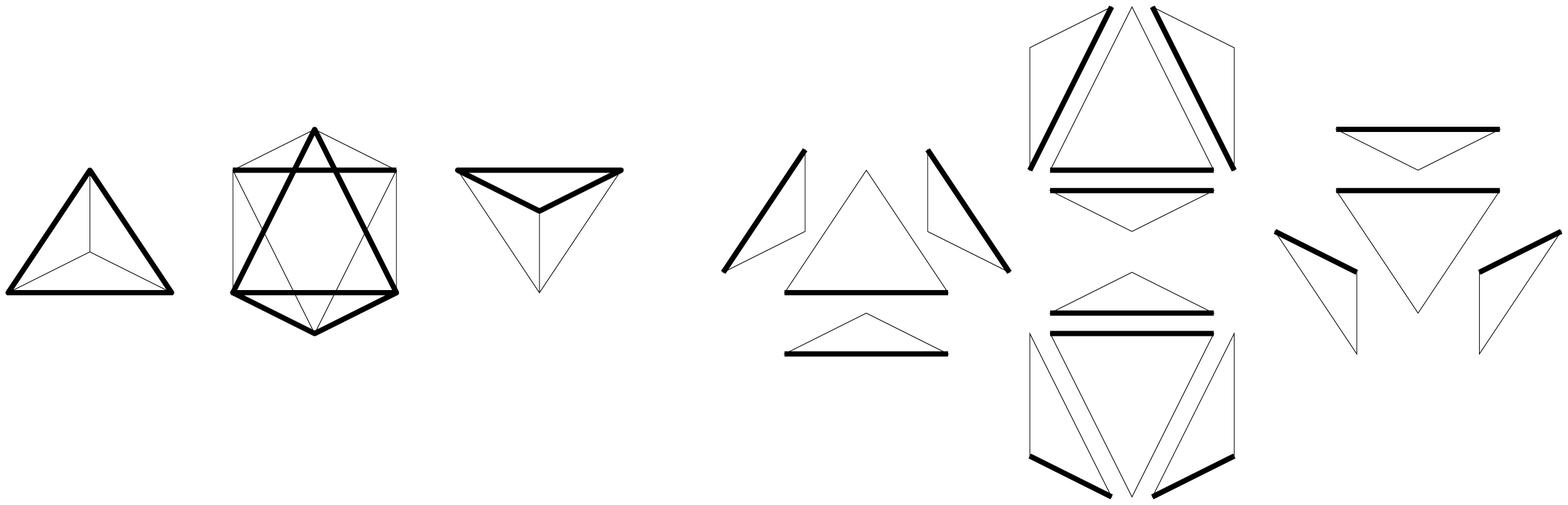}
\end{center}
\caption{\label{coreCDCnop.fig}$CC^{*}_{nop}D\quad (12)$.}
\end{figure}

\begin{figure}[ht]
\begin{center}
\includegraphics[scale = .4]{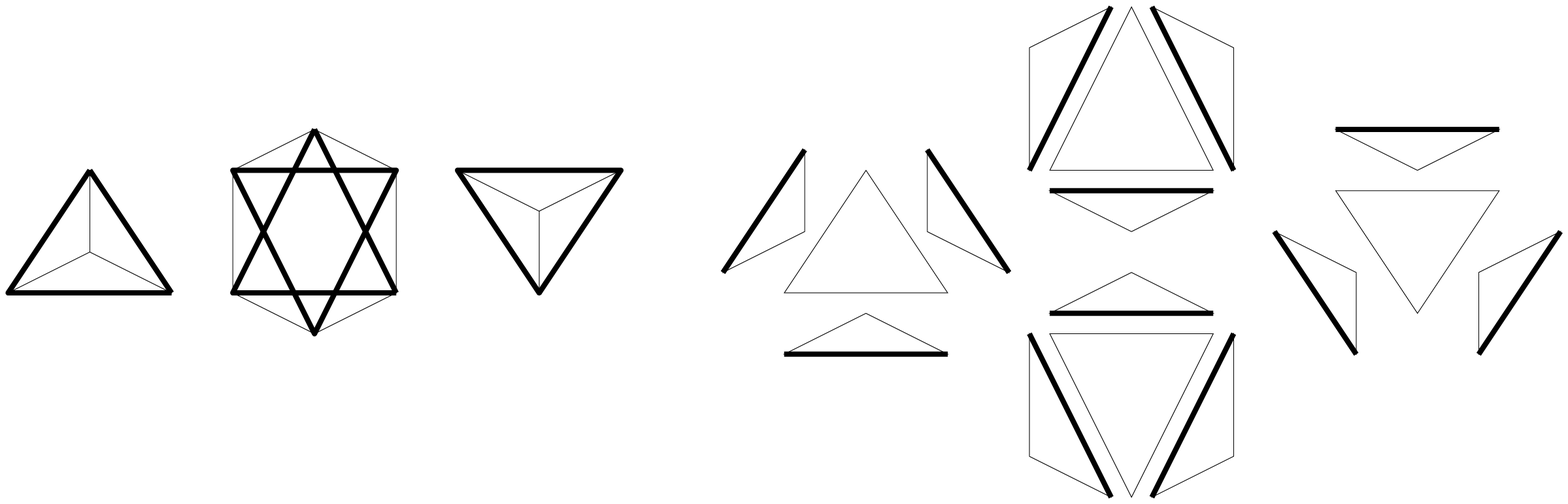}
\end{center}
\caption{\label{coreCDCop.fig}$CC^{*}_{op}D\quad (4)$.}
\end{figure}

\begin{proposition}
The subvariety $Z_{sp}$ consists of $66$ isolated points and $4$
subvarieties of positive dimension.  Modulo the action of the
symmetric group and the duality exchanging lines and $3$-spaces, there
are five types of points of $Z_{sp}$.  In terms of $\Gamma $, these
types appear in Figures~\ref{coreDDE.fig}--\ref{coreCDCop.fig}.
Figure~\ref{coreCDCop.fig} represents a point in the
positive-dimensional locus.  \footnote{The labels 
of the figures refer to certain divisors in $Z$, cf. \cite{tet.coho}.
The numbers in parentheses indicate how many of each type of component
appear, without modding out by the action of $S_{4}$.  The notation
$(a+a)$ indicates that there are $2a$ components of this type; we have
only depicted one of each dual pair.}
\end{proposition}

\begin{proof}
Let $\tilde{p}\in \tilde{U}$ be minimally split, and let $p$ be is its
image in $U$.  Then $n_{k} (\tilde{p}) = 2$ for $k\leq 3$, which
implies $\operatorname{Card}\{\pi _{I} (\tilde{p}) \}=6$.  Up to symmetry a
minimally split point must have its subspaces partitioned as follows:

\begin{enumerate}
\item The lines must collapse together as $(3,1)$ or $(2,2)$.  (The
notation $(p,q)$ means that the two distinct lines are the image of
$p$ and $q$ lines, where $p+q=4$.) 
\item The $2$-planes must collapse together as $(5,1)$, $(4,2)$, or
$(3,3)$. 
\item The $3$-planes, like the lines, must collapse together as
$(3,1)$ or $(2,2)$.  
\end{enumerate}

With these facts and Figures~\ref{core.fig} and \ref{vert-rels.fig} in
hand, computing $Z_{sp}$ becomes a combinatorial exercise.  We leave
the pleasure of this computation to the reader.
\end{proof}

\begin{theorem}\label{mainthm}
The core $\core $, and thus $\tilde{X}$, is nonsingular.      
\end{theorem}

\begin{proof}
We apply the Jacobian condition for nonsingularity in an affine
neighborhood of each point of $\core _{sp}$.  Let
$z=\eta(\tilde{p})\in\core_{sp}$.  Fix an affine neighborhood of $z$
in $\prod_{\beta\in\cH}\proj_{\beta}$ as follows.  For each
$\beta\in\cH$, choose one $\alpha(\beta)\in\cE(\beta)$ with
$y_{\alpha(\beta),\beta}\neq 0$, and set this coordinate equal to $1$.
The remaining variables $\{y_{\alpha,\beta}\mid
\alpha\neq\alpha(\beta)\}$ form a system of local parameters at $z$ in
$\prod_{\beta}\proj_{\beta}$.

Let $\Omega ^{1}_{z}$ be the $\C $-vector space of differentials of $\core
$ at the point $z$.  Since $\core $ is a threefold, we have $\dim
\Omega ^{1}_{z}\geq 3$, with equality only if $\core $ is nonsingular at $z$.
Furthermore,
\[
\Omega ^{1}_{z} = \Bigl(\bigoplus_{\substack{(\alpha ,\beta) \\
\alpha \not = \alpha(\beta)}}\, \C \cdot  dy_{\alpha,\beta}\Bigr)/{J},
\]
where $J$ is the subspace generated by the differentials (evaluated at
$z$) of all functions vanishing on $\core$.  To study this quotient,
we will use the following combinatorial rules for computing with
differentials in $\Omega ^{1}_{z}$.  These follow immediately from the
equations in Lemma~\ref{lem:core-equations}; we omit the simple proof.

\begin{lemma}\label{rules1}
Suppose that $y_{\alpha_1,\beta_1}y_{\alpha_2,\beta_2} -
y_{\alpha_2,\beta_1}y_{\alpha_1,\beta_2}$
vanishes on $\core $. %\vspace{.03in} 
\begin{enumerate}
\item If $y_{\alpha_1,\beta_1} = y_{\alpha_2,\beta_1} =
y_{\alpha_1,\beta_2} = 0$ and $y_{\alpha_2,\beta_2} \not =
0$, then $dy_{\alpha_1,\beta_1}=0$.%\vspace{.07in}
\item If $y_{\alpha_1,\beta_1} = y_{\alpha_1,\beta_2} = 0$ and
$y_{\alpha_2,\beta_2} = y_{\alpha_2,\beta_1}\not = 0$,
then $dy_{\alpha_1,\beta_1}=dy_{\alpha_1,\beta_2}$.
\end{enumerate}
\end{lemma}

Using these rules, we add data for $\Omega ^{1}_{z}$ to $\Gamma $ as
in Figure~\ref{rules1.fig}.  The $0$ means that the differential of
the variable corresponding to the edge is $0$, and the two $da$\,s
indicate that the two differentials coincide in $\Omega ^1_{z}$.

\begin{figure}[ht]
\begin{center}
\psfrag{0}{$\scriptscriptstyle{0}$}
\psfrag{da}{$\scriptscriptstyle{da}$}
\includegraphics[scale = .9]{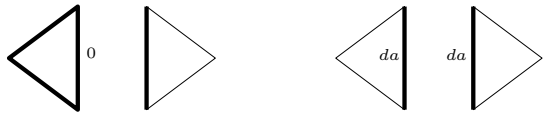}
\end{center}
\caption{\label{rules1.fig}Rules for differentials.}
\end{figure}

First we verify nonsingularity at the isolated points of $\core
_{sp}$.  Since the computations for the various points are all very
similar, we explain the case $DDE$ in detail and will leave the others
to the reader.  We fix the affine neighborhood of a point of type
$DDE$ by assigning the value $1$ to exactly one thin edge in each of
the $19$ components in Figure~\ref{coreDDE.fig}.  Since the
differentials of the linear polynomials are in $J$, the differential
corresponding to any thin edge is a linear combination of
differentials corresponding to bold edges.  Thus, $\Omega ^{1}_{z}$ is
generated by $dy_{\alpha,\beta}$, where $\alpha$ is a bold edge of the
component $\beta$.

Consider the hypersimplex connected components of
Figure~\ref{coreDDE.fig}.  Applying Lemma~\ref{rules1}, we find three
independent differentials $da$, $db$, and $dc$ in these components; the other differentials in these components are $0$.  Now consider
the other connected components of Figure~\ref{coreDDE.fig}.  Using
Lemma~\ref{rules1} we see that the remaining differentials are either
$0$ or are equal to $da $, $db$, or $dc$.  The result is summarized in
Figure~\ref{core-all-DDE.fig}.  Hence $\Omega ^{1}_{z}$ is
$3$-dimensional, and all the points of type $DDE$ are nonsingular
points of $\core $.  

\begin{figure}[ht]
\begin{center}
\psfrag{0}{$\scriptscriptstyle{0}$}
\psfrag{da}{$\scriptscriptstyle{da}$}
\psfrag{db}{$\scriptscriptstyle{db }$}
\psfrag{dc}{$\scriptscriptstyle{dc}$}
\includegraphics[scale = .6]{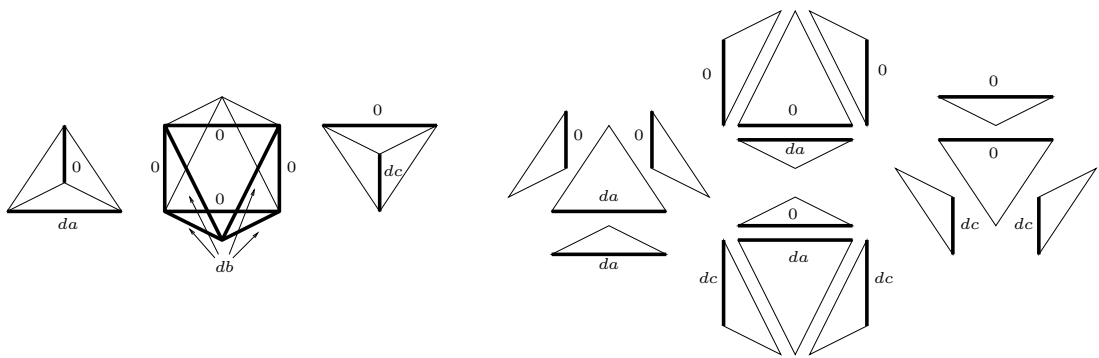}
\end{center}
\caption{\label{core-all-DDE.fig}All differentials for $DDE$.}
\end{figure}

Finally consider the family of special points $CC^{*}_{op}D$ in
Figure~\ref{coreCDCop.fig}.  In contrast to the isolated case, to
verify nonsingularity we have to use the cubic and quartic polynomials
of Lemma~\ref{lem:core-equations}.  Let $z\in Z$ be a point in a
subvariety of type $CC^{*}_{op}D$.  As before we construct an affine
neighborhood of $z$ choosing a thin edge in each connected component
of Figure~\ref{coreCDCop.fig} and setting
it to $1$.  At the point $z$, the linear relations imply that all of
the thin edges \emph{except those in the four thin triangles} will also
have value one.  To complete the graph $\Gamma $ to represent the point
$z$, we apply the quadratic relations to find $u,v\in \C$ such that
the values are as in Figure~\ref{ccd-1st.fig} (up to the choice of the
edges with values $1$).  Hence, a priori, this is a $2$-dimensional
component of $\core _{sp}$.

\begin{figure}[ht]
\begin{center}
\psfrag{2}{$\scriptscriptstyle{1}$}
\psfrag{x-1}{$\scriptscriptstyle{u-\frac{1}{2}}$}
\psfrag{-x-1}{$\scriptscriptstyle{-u-\frac{1}{2}}$}
\psfrag{y-1}{$\scriptscriptstyle{v-\frac{1}{2}}$}
\psfrag{-y-1}{$\scriptscriptstyle{-v-\frac{1}{2}}$}
\includegraphics[scale = .6]{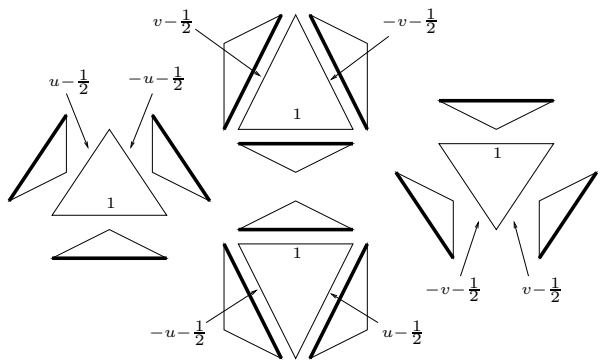}
\end{center}
\caption{\label{ccd-1st.fig}A typical point in $CC^{*}_{op}D$.}
\end{figure}

The quartic relations from Lemma~\ref{lem:core-equations},
however, imply that the two parameters $u$ and $v$ satisfy a linear
relation.  For the choice of parameters in Figure~\ref{ccd-1st.fig},
for example, this relation is 
\[
1\cdot(u-\frac{1}{2})\cdot 1\cdot 1 = 1\cdot (v-\frac{1}{2})\cdot
1\cdot 1, \quad \hbox{or}\quad u = v. 
\]
Thus this component of $\core _{sp}$ is in fact a curve.

We now complete the proof of the theorem.  As in the isolated case,
the differentials on all thin edges, except for those in the four thin
triangles, can be expressed as linear combinations of the differentials
on bold edges.  Moreover the differentials on edges of the thin triangles
can be expressed as linear combinations of $du$ and $dv$.  Thus, using
Lemma~\ref{rules1}, we can find $8$ differentials that span $\Omega
^{1}_{z}$ (Figure~\ref{core-ccd-2nd.fig}):
\[da , db , dc , da^*, db^*, dc^*, du, dv.\]
Note that the span of these is at most $5$-dimensional, because of the
relations $da + db + dc = da^* + db^* + dc^* = 0$ induced by the
differentials of the linear relations, and the relation $du=dv$
induced by the linear relation between $u$ and $v$.

\begin{figure}[ht]
\begin{center}
\psfrag{0}{$\scriptscriptstyle{0}$}
\psfrag{da}{$\scriptscriptstyle{da }$}
\psfrag{db}{$\scriptscriptstyle{db }$}
\psfrag{dc}{$\scriptscriptstyle{dc }$}
\psfrag{dA}{$\scriptscriptstyle{da^*}$}
\psfrag{dB}{$\scriptscriptstyle{db^* }$}
\psfrag{dC}{$\scriptscriptstyle{dc^* }$}
\includegraphics[scale = .6]{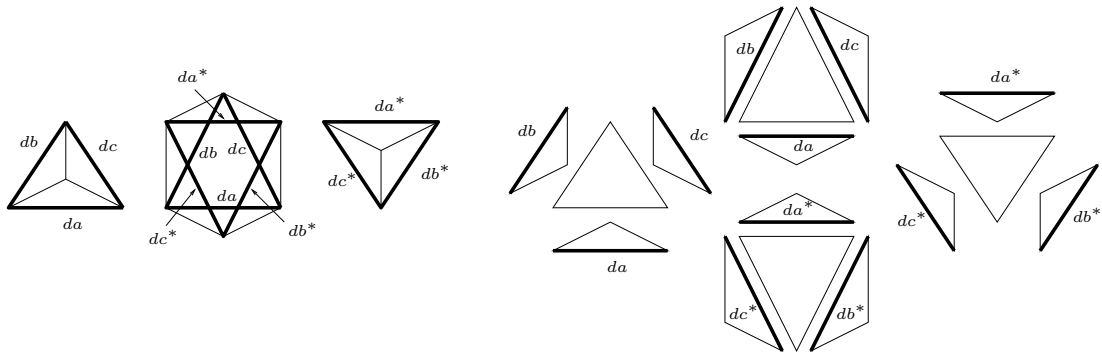}
\end{center}
\caption{\label{core-ccd-2nd.fig}Differentials for $CC^{*}_{op}D$.}
\end{figure}

To finish, we claim that the spans of $da ,
db , dc$ and $da^*, db^*, dc^*$ are each $1$-dimensional.  Indeed,
a quadric relation implies that the front face
of the octahedron in Figure~\ref{core-ccd-2nd.fig} has the same shape as
the corresponding face in Figure~\ref{ccd-1st.fig}, which implies 
\[
db  - (u-\frac{1}{2})\,dc  = 0.
\]
This relation, and a similar one involving $db^*$ and $dc^*$, shows that
$\Dim \Omega^{1}_{z} = 3$.  This completes the proof of the main theorem.
\end{proof}

\bibliographystyle{amsplain}
\bibliography{tetra}
\end{document}